\documentclass[11pt]{article}
\usepackage{latexsym,amssymb, color}
\usepackage{psfrag}
\usepackage{amsmath}

\newtheorem{Theorem}{\bf Theorem}[section]
\newtheorem{Lemma}{\bf Lemma}[section]
\newtheorem{Proposition}{\bf Proposition}[section]
\newtheorem{Corollary}{\bf Corollary}[section]
\newtheorem{Remark}{\bf Remark}[section]
\newtheorem{Example}{\bf Example}[section]
\newtheorem{Definition}{\bf Definition}[section]

\newenvironment{theorem}{\begin{Theorem}$\!\!\!$}{\end{Theorem}}
\newenvironment{lemma}{\begin{Lemma}$\!\!\!$}{\end{Lemma}}
\newenvironment{proposition}{\begin{Proposition}$\!\!\!$}{\end{Proposition}}
\newenvironment{corollary}{\begin{Corollary}$\!\!\!$}{\end{Corollary}}

\newenvironment{definition}{\begin{Definition}$\!\!\!$}{\end{Definition}}

\numberwithin{equation}{section}

\begin{document}

\title{Large time behavior of solutions of the heat equation\\ 
 with inverse square potential}
\author{Kazuhiro Ishige and Asato Mukai
}
\date{}
\maketitle
\begin{abstract}
Let $L:=-\Delta+V$ be a nonnegative Schr\"odinger operator on $L^2({\bf R}^N)$, 
where $N\ge 2$ and $V$ is a radially symmetric inverse square potential. 
In this paper we assume either $L$ is subcritical or null-critical 
and we establish a method for obtaining the precise description 
of the large time behavior of $e^{-tL}\varphi$, where $\varphi\in L^2({\bf R}^N,e^{|x|^2/4}\,dx)$. 
\end{abstract}
\vspace{40pt}
\noindent Addresses:

\smallskip
\noindent K. I.:  Mathematical Institute, Tohoku University,
Aoba, Sendai 980-8578, Japan.\\
\noindent 
E-mail: {\tt ishige@m.tohoku.ac.jp}\\

\smallskip
\noindent 
A. M.: Mathematical Institute, Tohoku University,
Aoba, Sendai 980-8578, Japan.\\
\noindent 
E-mail: {\tt asato.mukai.t7@dc.tohoku.ac.jp}\\
\vspace{20pt}
\vspace{3pt}
\newpage
\section{Introduction}
Let $L:=-\Delta+V$ be a nonnegative Schr\"odinger operator on $L^2({\bf R}^N)$, 
where $N\ge 2$ and $V$ is 
a radially symmetric inverse square potential, that is 
\begin{equation*}
\begin{split}
 & V(r)=\lambda_1r^{-2}+o(r^{-2+\theta})
\quad\mbox{as}\quad r\to 0,\\
 & V(r)=\lambda_2r^{-2}+o(r^{-2-\theta})
\quad\mbox{as}\quad r\to\infty,
\end{split}
\end{equation*}
for some $\lambda_1$, $\lambda_2\in[\lambda_*,\infty)$ with $\lambda_*:=-(N-2)^2/4$ and $\theta>0$. 
We are interested in the precise description 
of the large time behavior of $u=e^{-tL}\varphi$, which is a solution of 
\begin{equation}
\label{eq:1.1}
\left\{
\begin{array}{ll}
\partial_t u-\Delta u+V(|x|)u=0, & 
\quad x\in{\bf R}^N,\,\,t>0,\vspace{5pt}\\
u(x,0)=\varphi(x), & \quad x\in{\bf R}^N. 
\end{array}
\right.
\end{equation}
Nonnegative Schr\"odinger operators and their heat semigroups appear in various fields 
and have been studied intensively 
by many authors since the pioneering work due to Simon~\cite{S} 
(see e.g., \cite{Bar}, \cite{CK}, \cite{CUR}, \cite{DS}, 
\cite{IIY01}, \cite{IIY02}, \cite{IK02}--\cite{IKO}, 
\cite{LS}--\cite{MT2}, \cite{M0}, \cite{M}, 
\cite{P3}--\cite{Zhang} and references therein).  
See also the monographs of Davies \cite{Dav}, Grigor'yan \cite{Gri} and Ouhabaz \cite{Ouh}. 
The inverse square potential 
is a typical one appearing in the study of the Schr\"odinger operators 
and it arises in the linearized analysis for nonlinear diffusion equations
and in the asymptotic analysis for diffusion equations. 
\vspace{3pt}

Throughout this paper we assume the following condition on the potential~$V$: 
\begin{equation*}
(V)\qquad
\left\{
\begin{array}{ll}
({\rm i}) & \mbox{$V=V(r)\in C^1((0,\infty))$};\vspace{7pt}\\
({\rm ii}) &  \displaystyle{\lim_{r\to 0}r^{-\theta}\left|r^2V(r)-\lambda_1\right|=0},
 \quad
 \displaystyle{\lim_{r\to\infty}r^\theta\left|r^2V(r)-\lambda_2\right|=0},\vspace{3pt}\\
  & \mbox{for some $\lambda_1$, $\lambda_2\in[\lambda_*,\infty)$ with $\lambda_*:=-(N-2)^2/4$ and $\theta>0$};\vspace{7pt}\\
({\rm iii}) &  \displaystyle{\sup_{r\ge 1}\left|r^3V'(r)\right|<\infty}. 
\end{array}
\right.
\qquad
\end{equation*}
We say that $L:=-\Delta+V(|x|)$ is nonnegative on $L^2({\bf R}^N)$ if 
$$
\int_{{\bf R}^N}\left[|\nabla\phi|^2+V(|x|)\phi^2\right]\,dx\ge 0,
\qquad\phi\in C_0^\infty({\bf R}^N\setminus\{0\}). 
$$
When $L$ is nonnegative,  
we say that 
\begin{itemize}
  \item 
  $L$ is subcritical if, for any $W\in C_0({\bf R}^N)$, 
  $L-\epsilon W$ is nonnegative for all sufficiently small $\epsilon>0$:
  \item 
  $L$ is critical if $L$ is not subcritical. 
\end{itemize}
On the other hand, $L$ is said to be supercritical if $L$ is not nonnegative.  

Consider the ordinary differential equation
\begin{equation}
\tag{O}
U''+\frac{N-1}{r}U'-V(r)U=0\quad\mbox{in}\quad(0,\infty)
\end{equation}
under condition~(V). 
Equation~(O) has two linearly independent solutions $U$ (a {\it regular} solution) and $\tilde{U}$ (a~{\it singular} solution)
such that 
\begin{equation}
\label{eq:1.2}
U(r)\thicksim r^{A^+(\lambda_1)},
\qquad
\tilde{U}(r)\thicksim\left\{
\begin{array}{ll}
r^{A^-(\lambda_1)} & \mbox{if}\quad \lambda_1>\lambda_*,\vspace{3pt}\\
r^{-\frac{N-2}{2}}|\log r| & \mbox{if}\quad \lambda_1=\lambda_*,
\end{array}
\right.
\end{equation}
as $r\to +0$, where 
\begin{equation}
\label{eq:1.3}
A^\pm(\lambda):=\frac{-(N-2)\pm\sqrt{(N-2)^2+4\lambda}}{2}\quad\mbox{for}\quad\lambda\ge\lambda_*.
\end{equation}
In particular, $U\in L^2_{{\rm loc}}({\bf R}^N)$. 
Assume that $L$ is nonnegative on $L^2({\bf R}^N)$. 
Then it follows from \cite[Theorem~1.1]{IKO} that $U$ is positive in $(0,\infty)$  
and 
\begin{equation}
\label{eq:1.4}
U(r)\thicksim c_*v(r)\quad\mbox{as}\quad r\to\infty
\end{equation}
for some positive constant $c_*$, where
\begin{equation}
\label{eq:1.5}
v(r):=
\left\{
\begin{array}{ll}
r^{A^+(\lambda_2)} & \mbox{if $L$ is subcritical and $\lambda_2>\lambda_*$},\vspace{3pt}\\
r^{-\frac{N-2}{2}}\log r & \mbox{if $L$ is subcritical and $\lambda_2=\lambda_*$},\vspace{3pt}\\
r^{A^-(\lambda_2)} & \mbox{if $L$ is critical}.
\end{array}
\right.
\end{equation}
(See also \cite{M} for the case $\lambda_1=0$.) 
We often call $U$ a positive harmonic function for the operator~$L$.  
When $L$ is critical, 
following \cite{P2}, we say that 
$L$ is positive-critical if $U\in L^2({\bf R}^N)$
and that $L$ is null-critical if $U\not\in L^2({\bf R}^N)$. 
Generally, the behavior of the fundamental solution $p=p(x,y,t)$ corresponding to $e^{-tL}$ 
can be classified by whether $L$ is either subcritical, null-critical or positive-critical. 
Indeed, in the case of $\lambda_1=0$, 
by \cite[Theorem~1.2]{P2} we have: 
\begin{itemize}
  \item[(L1)] 
  If $L$ is subcritical, then $\displaystyle{\int_0^\infty p(x,y,t)\,dt<\infty}$ for $x$, $y\in{\bf R}^N$ with $x\not=y$;
  \item[(L2)] 
  If $L$ is null-critical, that is $A^-(\lambda_2)\ge -N/2$, 
  then   
  $$
  \int_0^\infty p(x,y,t)\,dt=\infty,\qquad
  \lim_{T\to\infty}\frac{1}{T}\int_0^T p(x,y,t)\,dt=0,
  $$ 
  for $x$, $y\in{\bf R}^N$ with $x\not=y$;
  \item[(L3)] 
  If $L$ is positive-critical, that is $A^-(\lambda_2)<-N/2$, then 
  $$
  \lim_{t\to\infty}p(x,y,t)=\frac{U(|x|)U(|y|)}{\,\,\,\,\,\|U\|_{L^2({\bf R}^N)}^2},
  \qquad x, y\in{\bf R}^N.
  $$
  \end{itemize}
See Corollary~\ref{Corollary:1.1} for (L1) and (L2) in the case of $\lambda_1\not=0$. 

On the other hand, under condition~$(V)$, 
the first author of this paper with Kabeya and Ouhabaz 
recently studied in \cite{IKO} the Gaussian estimate of the fundamental solution $p=p(x,y,t)$
in the subcritical case and in the critical case with $A^-(\lambda_2)>-N/2$. 
They proved that  
\begin{equation}
\label{eq:1.6}
0<p(x,y,t)
\le C\,t^{-\frac{N}{2}}
\frac{U(\min\{|x|,\sqrt{t}\})U(\min\{|y|,\sqrt{t}\})}{U(\sqrt{t})^2}
\exp\left(-\frac{|x-y|^2}{Ct}\right)
\end{equation}
holds for all $x$, $y\in{\bf R}^N$ and $t>0$, where $C$ is a positive constant (see \cite[Theorem~1.3]{IKO}). 
For related results, see e.g., \cite{Bar}, \cite{CUR}, \cite{GS}, \cite{LS}, \cite{MS1}, \cite{MT1}, \cite{MT2}, \cite{Zhang0}, 
\cite{Zhang} and references therein.  

The precise description of the large time behavior of $e^{-tL}\varphi$ with $\varphi\in L^2({\bf R}^N,e^{|x|^2/4}\,dx)$ 
has been studied in a series of papers~\cite{IK02}--\cite{IK05} only in the subcritical case 
with some additional restrictions 
such as $V\in C^1([0,\infty))$, $\lambda_2>\lambda_*$ and the sign of the potential. 
See also \cite{IK06}.
\vspace{3pt}

The purpose of this paper is to establish a method 
for obtaining the precise description of the large time behavior 
of $e^{-tL}\varphi$ with $\varphi\in L^2({\bf R}^N,e^{|x|^2/4}\,dx)$ 
in the subcritical case and in the null-critical case with $A^-(\lambda_2)>-N/2$, under condition~(V). 
In particular, we show that 
the solution~$u$ of \eqref{eq:1.1} behaves as a suitable multiple of 
$$
\left\{
\begin{array}{ll}
v_{{\rm reg}}(x,t) & \mbox{if $L$ is subcritical and $\lambda>\lambda_*$},\vspace{3pt}\\
(\log t)^{-1}v_{{\rm reg}}(x,t) & \mbox{if $L$ is subcritical and $\lambda=\lambda_*$},\vspace{3pt}\\
v_{{\rm sing}}(x,t) & \mbox{if $L$ is critical and $A^-(\lambda_2)>-N/2$},
\end{array}
\right.
$$
as $t\to\infty$ on all parabolic cones $\{x\in{\bf R}^N\,:\,R^{-1}t^{1/2}\le |x|\le Rt^{1/2}\}$ with $R>1$. 
(See Theorem~\ref{Theorem:1.4}.)
Here 
\begin{equation*}
\begin{split}
v_{{\rm reg}}(x,t):= & \,\,t^{-\frac{N+2A^+(\lambda_2)}{2}}
|x|^{A^+(\lambda_2)}\exp\left(-\frac{|x|^2}{4t}\right),\\
v_{{\rm sing}}(x,t):= & \,\,t^{-\frac{N+2A^-(\lambda_2)}{2}}
|x|^{A^-(\lambda_2)}\exp\left(-\frac{|x|^2}{4t}\right),
\end{split}
\end{equation*}
which are self-similar solutions of 
$$
\partial_t v=\Delta v-\lambda_2 |x|^{-2}v\quad\mbox{in}\quad{[\bf R}^N\setminus\{0\}]\times(0,\infty). 
$$
However, 
due to the fact that $v_{{\rm sing}}(t)\not\in H^1({\bf R}^N)$ for any $t>0$, 
the arguments in \cite{IK02}--\cite{IK05} are not applicable to the critical case. 
In this paper we study the large time behavior of the function 
$|x|^{-A}e^{-tL}\varphi$, instead of $e^{-tL}\varphi$, with
\begin{equation}
\label{eq:1.7}
\mbox{$A:=A^+(\lambda_2)$ if $L$ is subcritical},
\qquad
\mbox{$A:=A^-(\lambda_2)$ if $L$ is critical,}
\end{equation}
and overcome the difficulty arising from the fact 
that $v_{{\rm sing}}(t)\not\in H^1({\bf R}^N)$. 
As far as we know, this paper is the first one treating the precise large time behavior of $e^{-tL}\varphi$ 
in the critical case.  
\subsection{Radial solutions}
In this subsection 
we focus on radially symmetric solutions of \eqref{eq:1.1}. 
Divide the operator $L$ into the following three cases: 
$$
\begin{array}{ll}
& {\rm (S)}\,:\mbox{$L$ is subcritical and $\lambda>\lambda_*$};
\qquad\qquad
{\rm (S_*)}:\mbox{$L$ is subcritical and $\lambda=\lambda_*$};\vspace{3pt}\\
 & {\rm (C)}:\mbox{$L$ is critical and $A^-(\lambda_2)>-N/2$}.
\end{array}
$$ 
Set 
\begin{equation}
\label{eq:1.8}
d:=N+2A,\qquad
\rho_d(\xi):=\xi^{d-1}e^{\frac{\xi^2}{4}},\qquad
\psi_d(\xi):=c_de^{-\frac{\xi^2}{4}}. 
\end{equation}
Here $c_d$ is a positive constant such that $\|\psi_d\|_{L^2({\bf R}_+,\,\rho_d\,d\xi)}=1$, that is 
$c_d=[2^{d-1}\Gamma(d/2)]^{-1/2}$, where ${\bf R}_+:=(0,\infty)$ and $\Gamma$ is the Gamma function. 

Let $\varphi$ be radially symmetric and $\varphi\in L^2({\bf R}^N,e^{|x|^2/4}\,dx)$. 
Then $e^{-tL}\varphi$ is radially symmetric with respect to $x$ and 
set
$$
u(|x|,t)=\left[e^{-tL}\varphi\right](x),\qquad
v(|x|,t):=|x|^{-A}u(|x|,t),\qquad x\in{\bf R}^N,\,\,t>0. 
$$ 
Then $v$ satisfies the Cauchy problem for a {\it $d$-dimensional} parabolic equation
$$
\left\{
\begin{array}{ll}
\displaystyle{\partial_t v=\frac{1}{r^{d-1}}\partial_r\left(r^{d-1}\partial_r v\right)-V_{\lambda_2}(r)v}, & r\in(0,\infty),\,\,t>0,\vspace{5pt}\\
v(r,0)=r^{-A}\varphi(r), & r\in(0,\infty),
\end{array}
\right.
$$
where $V_{\lambda_2}(r):=V(r)-\lambda_2r^{-2}$. 
\vspace{5pt}

In the first and the second theorems 
we obtain the precise description of the large time behavior of the radially symmetric solutions of \eqref{eq:1.1} 
in either (S) or (C). 
\begin{theorem}
\label{Theorem:1.1}
Let $N\ge 2$ and assume condition~$(V)$. 
Let $L$ satisfy either {\rm (S)} or {\rm (C)}. 
Let $u=u(|x|,t)$ be a radially symmetric solution of \eqref{eq:1.1} 
such that $\varphi\in L^2({\bf R}^N,e^{|x|^2/4}\,dx)$. 
Define $w=w(\xi,s)$ by 
\begin{equation}
\label{eq:1.9}
w(\xi,s):=(1+t)^{\frac{d}{2}}r^{-A}u(r,t)\quad\mbox{with}\quad \xi=(1+t)^{-\frac{1}{2}}r\ge 0,\,\,\, s=\log(1+t)\ge 0.
\end{equation}
Then there exists a positive constant $C$ such that 
$$
\sup_{s>0}\|w(s)\|_{L^2({\bf R}_+,\rho_d\,d\xi)}\le C\|w(0)\|_{L^2({\bf R}_+,\rho_d\,d\xi)}.
$$
Furthermore, 
\begin{equation}
\label{eq:1.10}
\lim_{s\to\infty}w(\xi,s)=m(\varphi)\psi_d(\xi)
\quad\mbox{in}\quad
L^2({\bf R}_+,\rho_d\,d\xi)\,\cap\,C^2(K) 
\end{equation}
for any compact set $K$ in ${\bf R}^N\setminus\{0\}$, 
where
\begin{equation}
\label{eq:1.11}
m(\varphi):=\frac{c_d}{c_*}\int_0^\infty \varphi(r)U(r)r^{N-1}\,dr.
\end{equation}
In particular, if $m(\varphi)=0$, then 
\begin{equation}
\label{eq:1.12}
\|w(s)\|_{L^2({\bf R}_+,\rho_d\,d\xi)}+\|w(s)\|_{C^2(K)}=O(e^{-s})
\quad\mbox{as}\quad s\to\infty. 
\end{equation}
\end{theorem}
\begin{theorem}
\label{Theorem:1.2}
Assume the same conditions as in Theorem~{\rm\ref{Theorem:1.1}}. 
Set $u_*(r,t):=u(r,t)/U(r)$. 
\begin{itemize}
  \item[{\rm (a)}] 
  For any $j\in\{0,1,2\dots\}$, 
  $\partial_t^j u_*\in C([0,\infty)\times(0,\infty))$. 
  \item[{\rm (b)}]
  $\displaystyle{\lim_{t\to\infty}\,t^{\frac{d}{2}}u_*(0,t)=\frac{c_d}{c_*}m(\varphi)}$ 
  and 
  $\displaystyle{\lim_{t\to\infty}t^{\frac{d}{2}+1}(\partial_tu_*)(0,t)=-\frac{dc_d}{2c_*}m(\varphi)}$. 
  \item[{\rm (c)}]
  Let $T>0$ and $\epsilon$ be a sufficiently small positive constant. 
  Define 
  \begin{equation}
  \label{eq:1.13}
  G_d(r,t):=u_*(r,t)-\left[u_*(0,t)+(\partial_tu_*)(0,t)F_d(r)\right]\quad\mbox{for}\quad r\in[0,\infty),\,\,t>0,
  \end{equation}
  with
  $$
  U_d(s):=r^{-A}U(r),\quad
  F_d(r):=\int_0^r s^{1-d}[U_d(s)]^{-2}\left(\int_0^s \tau^{d-1}U_d(\tau)^2\,d\tau\right)\,ds.
  $$ 
  Then there exists a positive constant $C$ such that 
  \begin{equation}
  \label{eq:1.14}
  |(\partial_r^\ell G_d)(r,t)|\le Ct^{-\frac{d}{2}-2}r^{4-\ell}\|\varphi\|_{L^2({\bf R}^N,e^{|x|^2/4}\,dx)}
  \end{equation}
  for $\ell\in\{0,1,2\}$, $0\le r\le\epsilon(1+t)^{\frac{1}{2}}$ and $t\ge T$.  
  \end{itemize}
\end{theorem}
In case ($\mbox{S}_*$) we have: 
\begin{theorem}
\label{Theorem:1.3}
Let $N\ge 2$ and assume condition~$(V)$.  Let $L$ satisfy {\rm ($\mbox{S}_*$)}. 
Let $u=u(|x|,t)$ be a radially symmetric solution of \eqref{eq:1.1} 
such that $\varphi\in L^2({\bf R}^N,e^{|x|^2/4}\,dx)$. 
\begin{itemize}
  \item[{\rm (i)}] 
  Let $w$ be as in Theorem~{\rm\ref{Theorem:1.1}} and $K$ a compact set in ${\bf R}^N\setminus\{0\}$. 
  Then there exists a positive constant $C_1$ such that 
  $$
  \sup_{s>0}\,(1+s)\|w(s)\|_{L^2({\bf R}_+,\rho_2\,d\xi)}\le C_1\|w(0)\|_{L^2({\bf R}_+,\rho_2\,d\xi)}.
  $$
  Furthermore, 
  $$
  \lim_{s\to\infty}sw(\xi,s)=2m(\varphi)\psi_2(\xi)
  \quad\mbox{in}\quad
  L^2({\bf R}_+,\rho_2\,d\xi)\,\cap\,C^2(K),
  $$
  where $m(\varphi)$ is as in \eqref{eq:1.11}.
  \item[{\rm (ii)}] 
  Let $u_*$, $U_2$, $F_2$ and $G_2$ be as in Theorem~{\rm\ref{Theorem:1.2}} with $d=2$. 
  Then 
  \begin{equation*}
  \begin{split}
   & \partial_t^j u_*\in C([0,\infty)\times(0,\infty))
  \quad\mbox{for}\quad j\in\{0,1,2,\dots\},\\
   & \lim_{t\to\infty}t(\log t)^2 u_*(0,t)=2\sqrt{2}c_*^{-1}m(\varphi),\\
   & \lim_{t\to\infty}t^2(\log t)^2(\partial_t u_*)(0,t)=-2\sqrt{2}c_*^{-1}m(\varphi).
  \end{split}
  \end{equation*}
  Furthermore, 
  for any $T>0$ and any sufficiently small $\epsilon>0$, 
  there exists a positive constant $C_2$ such that 
  \begin{equation}
  \label{eq:1.15}
  |(\partial_r^\ell G_2)(r,t)|\le C_2t^{-3}[\log(2+t)]^{-2}r^{4-\ell}\|\varphi\|_{L^2({\bf R}^N,e^{|x|^2/4}\,dx)}
  \end{equation}
  for $\ell\in\{0,1,2\}$, $0\le r\le\epsilon(1+t)^{\frac{1}{2}}$ and $t\ge T$. 
\end{itemize}
\end{theorem}
%
\vspace{3pt}

The function $w$ defined by \eqref{eq:1.9} satisfies 
\begin{equation}
\label{eq:1.16}
\partial_s w+{\mathcal L}_d w+\tilde{V}(\xi,s)w=0
\quad\mbox{for}\quad \xi\in[0,\infty),\,\,s>0,
\end{equation}
where 
$$
{\mathcal L}_dw:=-\frac{1}{\rho_d(\xi)}\partial_\xi(\rho_d(\xi)\partial_\xi w)-\frac{d}{2}w,
\qquad
\tilde{V}(\xi,s):=e^sV_{\lambda_2}(e^{\frac{s}{2}}\xi).
$$ 
For the proof of Theorems~\ref{Theorem:1.1}--\ref{Theorem:1.3}, 
we regard the operator ${\mathcal L}_d$ as a {\it $d$-dimensional} elliptic operator with  
$$
\left\{
\begin{array}{ll}
d>2 & \mbox{in the case of (S)},\vspace{3pt}\\
d=2 & \mbox{in the case of $\lambda_2=\lambda_*$},\vspace{3pt}\\
0<d<2 & \mbox{in the case of (C) with $\lambda_2>\lambda_*$},
\end{array}
\right.
$$
and study the large time behavior of $w=w(\xi,s)$ by developing the arguments in a series of papers~\cite{I0}--\cite{IK05}. 
The function $\psi_d$ defined by \eqref{eq:1.8} is the first eigenfunction of the eigenvalue problem
\begin{equation}
\tag{E}
{\mathcal L}_d\phi=\mu\phi\quad\mbox{in}\quad{\bf R}_+,
\quad
\phi\in H^1({\bf R}_+,\rho_d(\xi)\,d\xi)
\end{equation}
and the corresponding eigenvalue is $0$ 
(see Lemma~\ref{Lemma:2.5}). 
We show that 
$w$ behaves like a suitable multiple of $\psi_d$ as $s\to\infty$. 
Furthermore, 
combining the radially symmetry of $u$ with the behavior of $w$, 
we prove Theorems~\ref{Theorem:1.1}--\ref{Theorem:1.3}. 

The eigenfunction $\psi_d$ corresponds to $v_{{\rm reg}}$ in the subcritical case 
and $v_{{\rm sing}}$ in the null-critical case, respectively.
In the null-critical case,  
$v_{reg}$ is transformed by \eqref{eq:1.9} into 
$$
e^{-\frac{A^+(\lambda_2)-A^-(\lambda_2)}{2}s}\tilde{\psi}_d
\quad
\mbox{with}\quad
\tilde{\psi}_d:=\xi^{A^+(\lambda_2)-A^-(\lambda_2)}e^{-\frac{\xi^2}{4}}.
$$
Here $\tilde{\psi}_d$ is the first eigenfunction of the eigenvalue problem 
$$
{\mathcal L}_d\phi=\mu\phi\quad\mbox{in}\quad{\bf R}_+,
\quad
\phi\in H^1_0({\bf R}_+,\rho_d(\xi)\,d\xi)
$$
and the corresponding eigenvalue is $[A^+(\lambda_2)-A^-(\lambda_2)]/2>0$. 
In the null-critical case with $\lambda_2>\lambda_*$, we see that 
$0<d<2$ and $H^1_0({\bf R}_+,\rho_d(\xi)\,d\xi)\not=H^1({\bf R}_+,\rho_d(\xi)\,d\xi)$. 
This justifies that 
the operator ${\mathcal L}_d$ has two positive eigenfunctions $\psi_d$ and $\tilde{\psi}_d$. 
\subsection{Nonradial solutions}
We discuss the large time behavior of solutions of \eqref{eq:1.1} 
without the radially symmetry of the solutions. 

Let $\Delta_{{\bf S}^{N-1}}$ be the Laplace-Beltrami operator on ${\bf S}^{N-1}$. 
Let $\{\omega_k\}_{k=0}^\infty$ be 
the eigenvalues of 
$$
-\Delta_{{\bf S}^{N-1}}Q=\omega Q\quad\mbox{on}\quad{\bf S}^{N-1},
\qquad
Q\in L^2({\bf S}^{N-1}). 
$$
Then $\omega_k=k(N+k-2)$ for $k=0,1,2,\dots$. 
Let
$\{Q_{k,i}\}_{i=1}^{\ell_k}$ and $\ell_k$ be 
the orthonormal system and the dimension of the eigenspace corresponding to $\omega_k$, respectively. 
Then, for any $\varphi\in L^2({\bf R}^N,e^{|x|^2/4}\,dx)$, 
we can find radially symmetric functions $\{\phi^{k,i}\}\subset L^2({\bf R}^N,e^{|x|^2/4}\,dx)$ such that 
$$
\varphi=\sum_{k=0}^\infty\sum_{i=1}^{\ell_k} \varphi^{k,i}\quad\mbox{in}\quad L^2({\bf R}^N,e^{|x|^2/4}\,dx),
\qquad
\varphi^{k,i}(x):=\phi^{k,i}(|x|)Q_{k,i}\left(\frac{x}{|x|}\right)
$$
(see \cite{I1} and \cite{IK02}). 
Define $L_k:=-\Delta+V_k(|x|)$ and $V_k(r):=V(r)+\omega_kr^{-2}$. 
Then 
\begin{equation}
\label{eq:1.17}
\begin{split}
\left[e^{-tL}\varphi^{k,i}\right](x) & =\left[e^{-tL_k}\phi^{k,i}\right](x)\,Q_{k,i}\left(\frac{x}{|x|}\right),\\
\left[e^{-tL}\varphi\right](x) & =\sum_{k=0}^\infty\sum_{i=1}^{\ell_k}\left[e^{-tL_k}\phi^{k,i}\right](x)\,Q_{k,i}\left(\frac{x}{|x|}\right)\\
 & \qquad\mbox{in $L^2({\bf R}^N)\cap L^\infty({\bf R}^N)$ for any $t>0$}.
\end{split}
\end{equation}
Therefore the behavior of $e^{-tL}\varphi$ is described 
by a series of the radially symmetric solutions $e^{-tL_k}\phi^{k,i}$. 
Furthermore, 
$V_k$ satisfies condition~(V) with $\lambda_1$ and $\lambda_2$ 
replaced by $\lambda_1+\omega_k$ and $\lambda_2+\omega_k$, respectively. 
In particular, $L_k$ is subcritical if $k\ge 1$. 
Therefore,
applying our results in Section~1.1, 
we can obtain the precise description of the large time behavior of $e^{-tL}\varphi$. 

As an application of the above argument, we obtain the following result. 
\begin{theorem}
\label{Theorem:1.4}
Let $N\ge 2$ and $\varphi\in L^2({\bf R}^N,e^{|x|^2/4}\,dx)$. 
Assume condition~$(V)$.  Let 
$$
M(\varphi):=\frac{1}{c_*\kappa}\int_{{\bf R}^N}\varphi(x)U(|x|)\,dx,
\quad
\kappa:=2^{N+2A}\pi^{\frac{N}{2}}\Gamma\left(\frac{N+2A}{2}\right)\biggr/\Gamma\left(\frac{N}{2}\right).
$$
\begin{itemize}
  \item[{\rm (a)}] 
  In cases {\rm (S)} and {\rm (C)}, 
  $$
  \lim_{t\to\infty}t^{\frac{N+A}{2}}[e^{-tL}\varphi](t^{\frac{1}{2}}y)=M(\varphi)|y|^Ae^{-\frac{|y|^2}{4}}
  $$
  in $L^2({\bf R}^N,e^{|y|^2/4}\,dy)$ and in $L^\infty(K)$ for any compact set $K\subset{\bf R}^N\setminus\{0\}$. 
  Furthermore, 
  $$
  \lim_{t\to\infty}t^{\frac{N+2A}{2}}\frac{[e^{-tL}\varphi](x)}{U(|x|)}=c_*^{-1}M(\varphi)
  $$
  uniformly on $B(0,R)$ for any $R>0$. 
  \item[{\rm (b)}] 
  In case {\rm ($\mbox{S}_*$)},
  $$
  \lim_{t\to\infty}t^{\frac{N+A}{2}}(\log t)\,[e^{-tL}\varphi](t^{\frac{1}{2}}y)=2M(\varphi)|y|^Ae^{-\frac{|y|^2}{4}}
  $$
  in $L^2({\bf R}^N,e^{|y|^2/4}\,dy)$ and in $L^\infty(K)$ for any compact set $K\subset{\bf R}^N\setminus\{0\}$. 
  Furthermore, 
  $$
  \lim_{t\to\infty}t^{\frac{N+2A}{2}}(\log t)^2\,\frac{[e^{-tL}\varphi](x)}{U(|x|)}=4c_*^{-1}M(\varphi)
  $$
  uniformly on $B(0,R)$ for any $R>0$. 
\end{itemize}
\end{theorem}
As a corollary of Theorem~\ref{Theorem:1.4}, we have: 
\begin{corollary}
\label{Corollary:1.1}
Let $N\ge 2$ and assume condition~{\rm (V)}. 
Let $x$, $y\in{\bf R}^N$. 
Then
\begin{equation*}
\begin{array}{ll}
\,\,\,\,\,
\displaystyle{\lim_{t\to\infty}t^{\frac{N+2A}{2}}\frac{p(x,y,t)}{U(|x|)U(|y|)}=(c_*^2\kappa)^{-1}}
 & \mbox{in cases {\rm (S)} and {\rm (C)}},\vspace{7pt}\\
\displaystyle{\lim_{t\to\infty}t(\log t)^2\frac{p(x,y,t)}{U(|x|)U(|y|)}=4(c_*^2\kappa)^{-1}}
 & \mbox{in case {\rm ($\mbox{S}_*$)}}.
\end{array}
\end{equation*}
\end{corollary}
Corollary~\ref{Corollary:1.1} implies the same conclusion as in (L1) and (L2). 
For related results, see e.g., \cite{CK}, \cite{MS1}, \cite{M0}, \cite{P0} and \cite{P2}. 

The above argument also enables us to obtain the higher order asymptotic expansions of $e^{-tL}\varphi$. 
Furthermore, similarly to \cite{I0}--\cite{IK05}, 
it is useful for the study the large time behavior of the {\it hot spots} of $e^{-tL}\varphi$. 
See a forthcoming paper.  
\vspace{5pt}

The rest of this paper is organized as follows. 
In Section~2 we formulate the definition of the solution of \eqref{eq:1.1} and prove some preliminary lemmas. 
In Section~3 we obtain a priori estimates of radially symmetric solutions of \eqref{eq:1.1} 
by using the comparison principle. 
In Section~4 we obtain the precise description of the large time behavior of radially symmetric solutions of \eqref{eq:1.1} 
and complete the proofs of Theorems~\ref{Theorem:1.1}--\ref{Theorem:1.3}. 
In Section~5, by the argument in Section~1.2 we apply Theorems~\ref{Theorem:1.1}--\ref{Theorem:1.3} 
to prove Theorem~\ref{Theorem:1.4} and Corollary~\ref{Corollary:1.1}. 
\section{Preliminaries}
We formulate the definition of the solution of \eqref{eq:1.1} 
and obtain some properties related to the operator $L$. 
For positive functions $f$ and $g$ defined in  $(0,R)$ for some $R>0$, 
we write 
$$
f(r)\thicksim g(r)\quad\mbox{as}\quad r\to 0
\quad\mbox{if}\quad
\lim_{r\to 0}\frac{f(r)}{g(r)}=1. 
$$ 
Similarly, for  positive functions $f$ and $g$ defined in  $(R,\infty)$ for some $R>0$, 
we write  
$$
f(r)\thicksim g(r)\quad\mbox{as}\quad r\to\infty
\quad\mbox{if}\quad
\lim_{r\to\infty}\frac{f(r)}{g(r)}=1. 
$$
By the letter $C$
we denote generic positive constants 
and they may have different values also within the same line. 
\subsection{Definition of the solution}
Assume condition~$(V)$ and let $L:=-\Delta+V$ be nonnegative. 
In this subsection we consider the Cauchy problem 
\begin{equation}
\tag{P}
\qquad
\left\{
\begin{array}{ll}
\partial_t u_*+L_*u_*=0 & \quad\mbox{in}\quad{\bf R}^N\times(0,\infty),
\vspace{5pt}\\
u_*(x,0)=\varphi_*(x) & \quad\mbox{in}\quad{\bf R}^N,
\vspace{3pt}
\end{array}
\right.
\end{equation}
where 
$$
L_*u_*:=-\frac{1}{\nu}\mbox{div}\,(\nu\nabla u_*),
\quad
\nu:=U^2\in L^1_{\rm loc}({\bf R}^N),
\quad
\varphi_*\in L^2({\bf R}^N,\,\nu\,dx).
$$
\begin{definition}
\label{Definition:2.1}
Let $\varphi_*\in L^2({\bf R}^N,\,\nu\,dx)$. 
We say that $u_*$ is a solution of $(P)$ if 
\begin{equation*}
\begin{split}
  & u_*\in C([0,\infty):L^2({\bf R}^N,\,\nu\,dx))\,\cap\,L^2((0,\infty):H^1({\bf R}^N,\,\nu\,dx)),\vspace{5pt}\\
  & \int_0^\infty\int_{{\bf R}^N}\left[-u_*\partial_th+\nabla u_*\nabla h\right]\nu\,dx=0
\quad\mbox{for any $h\in C^\infty_0({\bf R}^N\times(0,\infty))$},\\
  & \lim_{t\to +0}\|u_*(t)-\varphi_*\|_{L^2({\bf R}^N,\,\nu\,dx)}=0.
\end{split}
\end{equation*}
\end{definition}
Problem~(P) possesses a unique solution $u_*$ such that 
$$
\|u_*(t)\|_{L^2({\bf R}^N,\,\nu\,dx)}\le\|\varphi_*\|_{L^2({\bf R}^N,\,\nu\,dx)},
\qquad t>0,
$$
and we often denote by $e^{-tL_*}\varphi_*$ the unique solution~$u_*$. 
Since $U\in C^2({\bf R}^N\setminus\{0\})$ and $U>0$ in ${\bf R}^N\setminus\{0\}$, 
applying the parabolic regularity theorems (see e.g., \cite[Chapter~IV]{LSU}) to (P), we see that
\begin{equation}
\label{eq:2.1}
\partial_t^ju_*\in C^{2,1}([{\bf R}^N\setminus\{0\}]\times(0,\infty)),
\quad j=0,1,2,\dots.
\end{equation}
\begin{lemma}
\label{Lemma:2.1}
Assume condition~$(V)$ and that $L$ is nonnegative. 
Let $\varphi_*\in L^2({\bf R}^N,\,\nu\,dx)$ and $u_*:=e^{-tL_*}\varphi_*$. 
\begin{itemize}
  \item[{\rm (i)}] For any $j\in\{1,2,\dots\}$, 
  there exists $C>0$ such that 
  $$
  \|(\partial_t^ju_*)(t)\|_{L^2({\bf R}^N,\,\nu\,dx)}\le Ct^{-j}\|\varphi_*\|_{L^2({\bf R}^N,\,\nu\,dx)},
  \quad t>0.
  $$
  \item[{\rm (ii)}] 
  If $\varphi_*\in L^2({\bf R}^N,e^{|x|^2/4}\nu\,dx)$, then 
  \begin{equation*}
  \sup_{t>0}\|u_*(t)\|_{L^2({\bf R}^N,e^{|x|^2/4(1+t)}\nu\,dx)}\le\|\varphi_*\|_{L^2({\bf R}^N,e^{|x|^2/4}\nu\,dx)}.
  \end{equation*}
\end{itemize}
\end{lemma}
{\bf Proof.}
Assertion~(i) follows from the same argument as in the proof of \cite[Lemma~2.1]{I1}. 
We prove assertion~(ii). It follows that 
\begin{equation*}
\begin{split}
 & \int_0^t\int_{{\bf R}^N}(\partial_t u_*)u_*e^{\frac{|x|^2}{4(1+\tau)}}\nu\,dx\,d\tau\\
 & =\frac{1}{2}\int_{{\bf R}^N}u_*(x,\tau)^2e^{\frac{|x|^2}{4(1+\tau)}}\nu\,dx\biggr|_{\tau=0}^{\tau=t}
+\frac{1}{8}\int_0^t\int_{{\bf R}^N}u_*^2\frac{|x|^2}{(1+\tau)^2}e^{\frac{|x|^2}{4(1+\tau)}}\nu\,dx\,d\tau
\end{split}
\end{equation*}
and
\begin{equation*}
\begin{split}
 & \int_0^t\int_{{\bf R}^N}\nabla u_*\nabla\left[u_*e^{\frac{|x|^2}{4(1+\tau)}}\right]\nu\,dx\,d\tau\\
 & =\int_0^t\int_{{\bf R}^N}|\nabla u_*|^2 e^{\frac{|x|^2}{4(1+\tau)}}\nu\,dx\,d\tau
 +\int_0^t\int_{{\bf R}^N}u_*\nabla u_*\biggr[\frac{x}{2(1+\tau)}\biggr]e^{\frac{|x|^2}{4(1+\tau)}}\nu\,dx\,d\tau\\
 & \ge -\frac{1}{8}\int_0^t\int_{{\bf R}^N}u_*^2\frac{|x|^2}{(1+\tau)^2}e^{\frac{|x|^2}{4(1+\tau)}}\nu\,dx\,d\tau.
\end{split}
\end{equation*}
Then, by (P) we have 
\begin{equation*}
\int_{{\bf R}^N}u_*(x,t)^2e^{\frac{|x|^2}{4(1+\tau)}}\nu\,dx
\le\int_{{\bf R}^N}\varphi_*(x)^2e^{\frac{|x|^2}{4}}\nu\,dx
\end{equation*}
for $t>0$. Thus assertion~(ii) follows.  
(The proof of assertion~(ii) is somewhat formal, however it is justified 
by use of approximate solutions.)
$\Box$\vspace{3pt}
\newline
Furthermore, we have: 
\begin{lemma}
\label{Lemma:2.2}
Assume condition~{\rm (V)} and that $L$ is nonnegative. 
Let $u_*$ be a radially symmetric solution of {\rm (P)}. 
Then $\partial_t^ju_*$ is continuous in ${\bf R}^N\times(0,\infty)$, 
where $j\in\{0,1,2,\dots\}$.
\end{lemma}
{\bf Proof.}
Let $j\in\{0,1,2,\dots\}$ and set $v_j=\partial_t^ju_*$. 
By \eqref{eq:2.1} it suffices to prove the continuity of $v_j$ at $(0,t)\in{\bf R}^N\times(0,\infty)$.
Since  $v_j$ is radially symmetric, 
$v_j$ satisfies
\begin{equation}
\label{eq:2.2}
\begin{split}
\partial_t v_j
 & =\frac{1}{r^{N-1}\nu(r)}\partial_r(r^{N-1}\nu(r)\partial_r v_j)\\
 & =\frac{1}{r^{N+k-1}r^{-k}\nu(r)}\partial_r(r^{N+k-1}r^{-k}\nu(r)\partial_r v_j),
 \qquad r>0,\,\,t>0,
\end{split}
\end{equation}
for any $k\in{\bf R}$. 
Since $A^+(\lambda_1)>-N/2$, 
we can find $k\in\{1,2,\dots\}$ such that 
\begin{equation}
\label{eq:2.3}
-N-k<2A^+(\lambda_1)-k<N+k.
\end{equation}
Set 
$\tilde{v}_j({\bf x},t):=v_j(|{\bf x}|,t)$ and 
$\tilde{\nu}({\bf x}):=|{\bf x}|^{-k}\nu(|{\bf x}|)$
for ${\bf x}\in{\bf R}^{N+k}$ and $t>0$. 
By Definition~\ref{Definition:2.1}, Lemma~\ref{Lemma:2.1}~(i) and \eqref{eq:2.2} we see that 
$\tilde{v}_j$ satisfies
\begin{equation*}
\begin{split}
 & \partial_t\tilde{v}_j=\frac{1}{\tilde{\nu}}\mbox{div}_{N+k}\,(\tilde{\nu}\,\nabla_{N+k}\tilde{v}_j)
\quad\mbox{in}\quad{\bf R}^{N+k}\times(0,\infty),\\
 & \|\tilde{v}_j(t)\|_{L^2({\bf R}^{N+k},\,\tilde{\nu}\,d\bf{x})}
 =\|v_j(t)\|_{L^2({\bf R}^N,\,\nu\,dx)}\le Ct^{-j}\|\varphi_*\|_{L^2({\bf R}^N,\,\nu\,dx)}.
\end{split}
\end{equation*}
Furthermore, it follows from \eqref{eq:1.2} that 
$\tilde{\nu}({\bf x})\sim |{\bf x}|^{2A^+(\lambda_1)-k}$ as $|{\bf x}|\to 0$. 
This together with \eqref{eq:2.3} implies that $\tilde{\nu}$ is an $A_2$ weight in a neighborhood of ${\bf 0}\in{\bf R}^{N+k}$. 
By Lemma~\ref{Lemma:2.1}~(i), 
applying the regularity theorems for parabolic equations with $A_2$ weight 
(see e.g., \cite{CS} and \cite{I0}), we see that $\tilde{v}_j$ is continuous at $({\bf 0},t)\in{\bf R}^{N+k}\times(0,\infty)$. 
This means that $\partial_t^ju_*$ is continuous at $(0,t)\in{\bf R}^N\times(0,\infty)$. 
Thus Lemma~\ref{Lemma:2.2} follows.
$\Box$
\vspace{3pt}

We formulate the definition of the solution of \eqref{eq:1.1}. 
See also \cite{MT1} and \cite{MT2}. 
\begin{definition}
\label{Definition:2.2}
Let $u$ be a measurable function in ${\bf R}^N\times(0,\infty)$ and $\varphi\in L^2({\bf R}^N)$. 
Define 
$$
u_*(x,t):=\frac{u(x,t)}{U(|x|)},\qquad \varphi_*(x):=\frac{\varphi(x)}{U(|x|)}.
$$
Then we say that $u$ is a solution of \eqref{eq:1.1} if $u_*$ is a solution of {\rm (P)}. 
\end{definition}
In the case where $\lambda_1$, $\lambda_2>\lambda_*$, 
we can deduce from \eqref{eq:1.2} and \eqref{eq:1.3} that 
$U\in H^1({\bf R}^N)$ and that
a solution $u$ of \eqref{eq:1.1} satisfies 
$$
u\in C([0,\infty):L^2({\bf R}^N))\,\cap\,L^2((0,\infty):H^1({\bf R}^N)). 
$$
We remark that $\varphi\in L^2({\bf R}^N)$ if and only if $\varphi_*\in L^2({\bf R}^N,\,\nu\,dx)$. 
Furthermore, by \eqref{eq:1.6} we have the following lemma
(see also \cite[Theorem~1.2]{IIY01} and \cite[Theorem~1.1]{IIY02}). 
\begin{lemma}
\label{Lemma:2.3}
Let $u$ be a solution of \eqref{eq:1.1} under condition~{\rm (V)}.
Assume either $L$ is subcritical or $L$ is critical with $A^-(\lambda_2)>-N/2$. 
Then, for any $T>0$, there exists  $C>0$ such that 
\begin{equation}
\label{eq:2.4}
\frac{|u(x,t)|}{U(\min\{|x|,\sqrt{t}\})}\le Ct^{-\frac{N}{4}}U(\sqrt{t})^{-1}\|\varphi\|_{L^2({\bf R}^N)},
\quad x\in{\bf R}^N,\,\,t\ge T.
\end{equation}
\end{lemma}
{\bf Proof.}
It follows from \eqref{eq:1.6} that 
\begin{equation*}
\begin{split}
 & \frac{|u(x,t)|}{U(\min\{|x|,\sqrt{t}\})}
 \le\frac{1}{U(\min\{|x|,\sqrt{t}\})}\left(\int_{\{|y|\le\sqrt{t}\}}+\int_{\{|y|>\sqrt{t}\}}\right)p(x,y,t)|\varphi(y)|\,dy\\
 & \le Ct^{-\frac{N}{2}}U(\sqrt{t})^{-2}\int_{\{|y|\le\sqrt{t}\}}|\varphi(y)|U(|y|)\,dy
+Ct^{-\frac{N}{2}}U(\sqrt{t})^{-1}\int_{\{|y|>\sqrt{t}\}}e^{-\frac{|x-y|^2}{Ct}}|\varphi(y)|\,dy\\
 & \le Ct^{-\frac{N}{2}}U(\sqrt{t})^{-2}\|U\|_{L^2(\{|y|\le\sqrt{t}\})}\|\varphi\|_{L^2({\bf R}^N)}
 +Ct^{-\frac{N}{4}}U(\sqrt{t})^{-1}\|\varphi\|_{L^2({\bf R}^N)}
\end{split}
\end{equation*}
for $x\in{\bf R}^N$ and $t>0$. 
On the other hand, 
by \eqref{eq:1.4} and \eqref{eq:1.5} 
we have
$$
\|U\|_{L^2(\{|y|\le\sqrt{t}\})}\le Ct^{\frac{N}{4}}U(\sqrt{t})
$$
for $t\ge T$ (see also \eqref{eq:3.7}). 
These imply \eqref{eq:2.4} and Lemma~\ref{Lemma:2.3} follows. 
$\Box$
\subsection{Preliminary lemmas}
We prove a lemma on 
the decay of $U'$ as $r\to\infty$. 
\begin{lemma}
\label{Lemma:2.4}
Let $N\ge 2$. 
Assume condition~{\rm (V)} and that $L=-\Delta+V(|x|)$ is nonnegative. 
Let $U$ and $v$ be as in \eqref{eq:1.2} and \eqref{eq:1.5}, respectively. 
In cases {\rm (S)} and {\rm (C)} there exists $\delta>0$ such that
\begin{equation}
\label{eq:2.5}
[v(r)^{-1}U(r)]'=O(r^{-1-\delta})\quad\mbox{as}\quad r\to\infty.
\end{equation}
\end{lemma}
{\bf Proof.}
Let $V_{\lambda_2}(r):=V(r)-\lambda_2r^{-2}$. Set
$$
v^+(r):=
\left\{
\begin{array}{ll}
r^{-\frac{N-2}{2}}\log r & \mbox{if $L$ is subcritical and $\lambda=\lambda_*$},\vspace{3pt}\\
r^{A^+(\lambda_2)} & \mbox{otherwise},
\end{array}
\right.
\qquad
v^-(r):=r^{A^-(\lambda_2)}.
$$
It follows from \eqref{eq:1.4} and (V)~(ii) that 
\begin{equation}
\label{eq:2.6}
\begin{split}
 & \tau^{N-1}v^-(\tau)V_{\lambda_2}(\tau)U(\tau)=O(\tau^{N-3-\theta+A^-(\lambda_2)}v(\tau))\\
 & =\left\{
\begin{array}{ll}
O(\tau^{-1-\theta}) & \mbox{if $L$ is subcritical and $\lambda_2>\lambda_*$},\vspace{3pt}\\
O(\tau^{-1-\theta-\sqrt{Q}}) & \mbox{if $L$ is critical and $\lambda_2>\lambda_*$},\vspace{3pt}\\
O(\tau^{-1-\theta}) & \mbox{if $L$ is critical and $\lambda_2=\lambda_*$},
\end{array}
\right.
\end{split}
\end{equation}
as $\tau\to\infty$, where $Q=(N-2)^2+4\lambda_2$. 
Then the function
$$
G(r):=v^-(r)
\int_1^r s^{1-N}[v^-(s)]^{-2}\left(\int_s^\infty\tau^{N-1}v^-(\tau)V_{\lambda_2}(\tau)U(\tau)\,d\tau\right)\,ds
$$
can be defined for any $r>0$ and satisfies 
\begin{eqnarray}
\nonumber
 & & G''(r)+\frac{N-1}{r}G'(r)-\lambda_2 r^{-2}G(r)=V_{\lambda_2}(r)U(r)\quad\mbox{in}\quad(0,\infty),\\
\label{eq:2.7}
 & & G(r)=o(v^+(r))\quad\mbox{as}\quad r\to\infty. 
\end{eqnarray}
Since 
$$
U''(r)+\frac{N-1}{r}U'(r)-\lambda_2r^{-2}U(r)=V_{\lambda_2}(r)U(r)\quad\mbox{in}\quad(0,\infty),
$$
the function $\tilde{v}(r):=U(r)-G(r)$ satisfies
\begin{equation}
\label{eq:2.8}
\tilde{v}''(r)+\frac{N-1}{r}\tilde{v}'(r)-\lambda_2 r^{-2}\tilde{v}(r)=0\quad\mbox{in}\quad(0,\infty). 
\end{equation}
On the other hand, 
$v^\pm$ satisfy \eqref{eq:2.8} and are linearly independent. 
Therefore, applying the standard theory for ordinary differential equations, 
we can find $a$, $b\in{\bf R}$ such that 
$\tilde{v}(r)=av^+(r)+bv^-(r)$ in $(0,\infty)$, that is 
\begin{equation}
\label{eq:2.9}
U(r)=av^+(r)+bv^-(r)+G(r)\quad\mbox{in}\quad(0,\infty). 
\end{equation}

Assume that $L$ is subcritical.  
By \eqref{eq:1.4}, \eqref{eq:2.7} and \eqref{eq:2.9}
we have 
\begin{equation*}
\begin{split}
 & v(r)^{-1}U(r)\\
 & =c_*+r^{-\sqrt{Q}}
\left[b+\int_1^r s^{1-N}[v^-(s)]^{-2}\left(\int_s^\infty\tau^{N-1}v^-(\tau)V_{\lambda_2}(\tau)U(\tau)\,d\tau\right)\,ds\right].
\end{split}
\end{equation*}
Since $Q=(N-2)^2+4\lambda_2>0$, 
by \eqref{eq:2.6} 
we can find $\delta'>0$ such that 
\begin{equation*}
\begin{split}
 & \left[v(r)^{-1}U(r)\right]'\\
 & =-\sqrt{Q}r^{-\sqrt{Q}-1}\left[b+\int_1^r s^{1-N}[v^-(s)]^{-2}\left(\int_s^\infty\tau^{N-1}v^-(\tau)V_{\lambda_2}(\tau)U(\tau)\,d\tau\right)\,ds\right]\\
 & \qquad\quad
+r^{-\sqrt{Q}}r^{1-N}[v^-(r)]^{-2}\int_r^\infty\tau^{N-1}v^-(\tau)V_{\lambda_2}(\tau)U(\tau)\,d\tau\\
 & =O(r^{-1-\delta'})\quad\mbox{as}\quad r\to\infty. 
\end{split}
\end{equation*}
This implies \eqref{eq:2.5} in the subcritical case. 

Next we assume that $L$ is critical. 
By \eqref{eq:1.4} and \eqref{eq:2.7} we see that $a=0$ and 
\begin{equation*}
\begin{split}
v(r)^{-1}U(r)
=b+\int_1^r s^{1-N}[v^-(s)]^{-2}\left(\int_s^\infty\tau^{N-1}v^-(\tau)V_{\lambda_2}(\tau)U(\tau)\,d\tau\right)\,ds.
\end{split}
\end{equation*}
This together with \eqref{eq:2.6} implies that 
$$
\left[v(r)^{-1}U(r)\right]'
=r^{1-N}[v^-(r)]^{-2}\left(\int_r^\infty\tau^{N-1}v^-(\tau)V_{\lambda_2}(\tau)U(\tau)\,d\tau\right)\,ds
=O(r^{-1-\theta})
$$
as $r\to\infty$, and \eqref{eq:2.5} holds with $\delta=\theta$. 
Thus Lemma~\ref{Lemma:2.4} follows. 
$\Box$\vspace{5pt}

At the end of this section we state the following lemma on eigenvalue problem~(E). 
\begin{lemma}
\label{Lemma:2.5}
Let $\{\mu_i\}_{i=0}^\infty$ be the eigenvalues of {\rm (E)} such that $\mu_0\le\mu_1\le\mu_2\le\dots$. 
Then, for any $i\in\{0,1,2,\dots\}$, 
$\mu_i=i$ and $\mu_i$ is simple. 
Furthermore, $\psi_d$ given in \eqref{eq:1.8} 
is the first eigenfunction of {\rm (E)}.
\end{lemma}
{\bf Proof.}
We leave the proof to the reader since 
it is proved by the same argument as in \cite[Lemma~2.1]{MNY}.
$\Box$
\section{A priori estimates of radial solutions}
Let $T>0$ and $\epsilon>0$. Define
$$
D_\epsilon(T):=\left\{(x,t)\in{\bf R}^N\times(T,\infty)\,:\,|x|<\epsilon t^{\frac{1}{2}}\right\}.
$$
In this section we prove the following proposition. 
\begin{proposition}
\label{Proposition:3.1}
Assume condition~$(V)$. Let $L$ satisfy either {\rm (S)}, {\rm $(\mbox{S}_*)$} or {\rm (C)}. 
Let $u_*=u_*(|x|,t)$ be a radially symmetric solution of {\rm (P)} such that 
$\|\varphi_*\|_{L^2({\bf R}^N,\,\nu\,dx)}=1$. 
Assume that
\begin{equation}
\label{eq:3.1}
\sup_{t>0}\,t^D[\log (2+t)]^{D'}\|u_*(t)\|_{L^2({\bf R}^N,\,\nu\,dx)}<\infty
\,\,\,\,\mbox{for some $D\ge 0$ and $D'\ge 0$}.
\end{equation}
Let $j\in\{0,1,2,\dots\}$. Then the following holds 
for any $T>0$ and any sufficiently small $\epsilon>0$. 
\begin{itemize}
  \item[{\rm (i)}] 
  There exists $C_1>0$ such that
  $$
  |(\partial_t^ju_*)(|x|,t)|\le C_1\Gamma_{D,D',j}(t)
  $$
  for $(x,t)\in D_\epsilon(T)$, where 
  \begin{equation}
  \label{eq:3.2}
  \Gamma_{D,D',j}(t):=
  \left\{
  \begin{array}{ll}
  t^{-D-\frac{d}{4}-j}[\log(2+t)]^{-D'} & \mbox{in the case of {\rm (S)}},\vspace{3pt}\\
  t^{-D-\frac{d}{4}-j}[\log(2+t)]^{-D'-1} & \mbox{in the case of {\rm ($\mbox{S}_*$)}},\vspace{3pt}\\
  t^{-D-\frac{d}{4}-j}[\log(2+t)]^{-D'} & \mbox{in the case of {\rm (C)}}.\vspace{3pt}
  \end{array}
  \right.
  \end{equation}
  \item[{\rm (ii)}] 
  Let 
  $$
  F_N^j(r,t):=\int_0^r s^{1-N}[\nu(s)]^{-1}\left(\int_0^s \tau^{N-1}\nu(\tau)(\partial_t^{j+1} u_*)(\tau,t)\,d\tau\right)\,ds.
  $$
  Then
  $$
  (\partial_t^ju_*)(|x|,t)=(\partial_t^ju_*)(0,t)+F_N^j(|x|,t)\quad\mbox{in}\quad{\bf R}^N\times(0,\infty). 
  $$
  Furthermore, there exists $C_2>0$ such that  
  $$
  |F_N^j(|x|,t)|\le C_2\Gamma_{D,D',j+1}(t)|x|^2,
  \qquad
  |(\partial_r F_N^j)(|x|,t)|\le C\Gamma_{D,D',j+1}(t)|x|,
  $$
  for $(x,t)\in D_\epsilon(T)$. 
\end{itemize}
\end{proposition}
For the proof, 
we construct supersolutions of problem~(P) in $D_\epsilon(T)$. 
\begin{lemma}
\label{Lemma:3.1}
Assume condition~$(V)$. 
Let $\gamma_1\ge 0$ and $\gamma_2\ge 0$. 
Set
$$
\zeta(t):=t^{-\gamma_1}[\log(2+t)]^{-\gamma_2}.
$$
Then, for any $T>0$ and any sufficiently small $\epsilon>0$,  
there exists a function $W_*=W_*(x,t)$ such that 
\begin{eqnarray}
\label{eq:3.3}
 & & \partial_t W_*+L_*W_*\ge 0
\qquad\quad\,\mbox{in}\quad{\bf R}^N\times(0,\infty),\\
\label{eq:3.4}
 & & \zeta(t)\le W_*(x,t)\le 2\zeta(t)
\quad\mbox{in}\quad D_\epsilon(T). 
\end{eqnarray}
\end{lemma}
{\bf Proof.}
Let $T>0$ and $\epsilon>0$. 
Let $\kappa$ be a positive constant such that 
\begin{equation}
\label{eq:3.5}
|\zeta'(t)|\le\kappa t^{-1}\zeta(t),\qquad t>0.
\end{equation}
Let 
$$
F(x):=\int_0^{|x|}s^{1-N}[\nu(s)]^{-1}\left(\int_0^s \tau^{N-1}\nu(\tau)\,d\tau\right)\,ds,
$$
which satisfies $-L_*F=1$ in ${\bf R}^N$. 
Set 
\begin{equation}
\label{eq:3.6}
W_*(x,t):=2\zeta(t)\left[1-\kappa t^{-1}F(x)\right]. 
\end{equation}
Since $\zeta$ is monotone decreasing, 
by \eqref{eq:3.5} we have 
\begin{equation*}
\begin{split}
 & \partial_t W_*+L_*W_*\\
 & \ge 2\zeta'(t)\left[1-\kappa t^{-1}F(x)\right]
+2\kappa\zeta(t) t^{-2}F(x)
+2\kappa t^{-1}\zeta(t)\\
& \ge 2\zeta'(t)+2\kappa t^{-1}\zeta(t)\ge 0
\quad\mbox{in}\quad{\bf R}^N\times(0,\infty). 
\end{split}
\end{equation*}
This implies \eqref{eq:3.3}. 
On the other hand, 
by \eqref{eq:1.2}, \eqref{eq:1.4} and \eqref{eq:1.5}
we have
\begin{equation*}
\begin{split}
 & \int_0^s \tau^{N-1}\nu(\tau)\,d\tau\le Cs^{2A^+(\lambda_1)+N}\quad\mbox{for}\quad 0<s\le 1,\\
 & \int_0^s \tau^{N-1}\nu(\tau)\,d\tau\le
\left\{
\begin{array}{ll}
s^{2A^+(\lambda_2)+N} & \mbox{in the case of (S)},\vspace{3pt}\\
s^2[\log(2+s)]^2 & \mbox{in the case of ($\mbox{S}_*$)},\vspace{3pt}\\
s^{2A^-(\lambda_2)+N} & \mbox{in the cases of (C)},
\end{array}
\right.
\quad\mbox{for}\quad s>1.
\end{split}
\end{equation*}
These imply that 
\begin{equation}
\label{eq:3.7}
\int_0^s \tau^{N-1}\nu(\tau)\,d\tau\le Cs^N\nu(s),\qquad s\ge 0. 
\end{equation}
Then it follows that $0\le F(x)\le C|x|^2$ for $x\in{\bf R}^N$. 
Taking a sufficiently small $\epsilon>0$ if necessary, 
we obtain 
$$
0\le\kappa t^{-1}F(x)\le C\epsilon^2\kappa\le\frac{1}{2},
\quad
(x,t)\in D_\epsilon(T). 
$$
This together with \eqref{eq:3.6} implies \eqref{eq:3.4}. 
Thus Lemma~\ref{Lemma:3.1} follows. 
$\Box$
\vspace{5pt}

Applying the same argument as in \cite[Lemma~3.2]{IK05}, we have: 
\begin{lemma}
\label{Lemma:3.2}
Assume the same conditions as in Proposition~{\rm\ref{Proposition:3.1}}. 
Furthermore, assume \eqref{eq:3.1} for some $D\ge 0$ and $D'\ge 0$. 
Let $T>0$ and let $\epsilon$ be a sufficiently small positive constant. 
Then, for any $j\in\{0,1,2,\dots\}$, 
there exists $C>0$ such that
\begin{equation}
\label{eq:3.8}
|(\partial_t^ju_*)(|x|,t)|\le C\Gamma_{D,D',j}(t)
\quad\mbox{in}\quad D_\epsilon(T).
\end{equation}
\end{lemma}
{\bf Proof.}
Let $j\in\{0,1,2,\dots\}$. Set 
$v_j:=\partial_t^j u_*$ and $u_j:=U(|x|)v_j(x,t)$. 
Since 
$$
v_j(\cdot,t)=\partial_t^j\left[e^{-(t/2)L_*}u_*(t/2)\right],\quad t>0,
$$
Lemma~\ref{Lemma:2.1} together with \eqref{eq:3.1} implies that 
$$
\sup_{t>0}\,t^{D+j}[\log (2+t)]^{D'}\|u_j(t)\|_{L^2({\bf R}^N)}
=\sup_{t>0}\,t^{D+j}[\log (2+t)]^{D'}\|v_j(t)\|_{L^2({\bf R}^N,\,\nu\,dx)}<\infty. 
$$
Let $T>0$ and let $\epsilon$ be a sufficiently small positive constant. 
Since $u_j$ satisfies 
$$
\partial_t u_j=\Delta u_j-V(|x|)u_j\quad\mbox{in}\quad{\bf R}^N\times(0,\infty), 
$$
by Lemma~\ref{Lemma:2.3} we have 
$$
|u_j(|x|,t)|\le Ct^{-\frac{N}{4}}\|u_j(t/2)\|_{L^2({\bf R}^N)}
\le Ct^{-D-\frac{N}{4}-j}[\log(2+t)]^{-D'}
$$
for all $x\in{\bf R}^N$ and $t>T$ with $|x|\ge\epsilon(1+t)^{1/2}$. 
This together with \eqref{eq:1.4}, \eqref{eq:1.5}, \eqref{eq:1.7}, \eqref{eq:1.8} and \eqref{eq:3.2}
implies that 
\begin{equation}
\label{eq:3.9}
|v_j(|x|,t)|\le C\frac{u_j(|x|,t)}{U(\epsilon(1+t)^{\frac{1}{2}})}
\le C\Gamma_{D,D',j}(t)
\end{equation}
for all $(x,t)\in{\bf R}^N\times[T,\infty)$ with $|x|=\epsilon(1+t)^{\frac{1}{2}}$. 
On the other hand, 
it follows from Lemma~\ref{Lemma:2.2} that 
\begin{equation}
\label{eq:3.10}
|v_j(|x|,T)|\le C\quad\mbox{for $x\in{\bf R}^N$ with $|x|\le\epsilon(1+T)^{\frac{1}{2}}$}. 
\end{equation}
Let $W_*$ be as in Lemma~\ref{Lemma:3.1} with $\zeta$ replaced by $\Gamma_{D,D',j}$. 
Then, by Lemma~\ref{Lemma:3.1}, \eqref{eq:3.9} and \eqref{eq:3.10} 
we apply the comparison principle to obtain 
$$
|v_j(|x|,t)|\le CW_*(x,t)\le 2C\Gamma_{D,D',j}(t)\quad\mbox{in}\quad D_\epsilon(T). 
$$
This implies \eqref{eq:3.8}, and the proof is complete. 
$\Box$\vspace{5pt}

Now we are ready to complete the proof of Proposition~\ref{Proposition:3.1}.
\vspace{3pt}
\newline
{\bf Proof of Proposition~\ref{Proposition:3.1}.}
By Lemma~\ref{Lemma:3.2} 
it suffices to prove assertion~(ii). 
Let $T>0$ and let $\epsilon$ be a sufficiently small positive constant. 
By \eqref{eq:3.7} and \eqref{eq:3.8}
we obtain  
\begin{equation*}
\begin{split}
 & |F_N^j(|x|,t)|\le C\Gamma_{D,D',j+1}(t)
\int_0^{|x|} s^{1-N}[\nu(s)]^{-1}\left(\int_0^s \tau^{N-1}\nu(\tau)\,d\tau\right)\,ds\\
 & \qquad\qquad\,\,\,\,
\le C\Gamma_{D,D',j+1}(t)|x|^2,\\
 & |(\partial_r F_N^j)(|x|,t)|\le C\Gamma_{D,D',j+1}(t)|x|,
\end{split}
\end{equation*}
for $(x,t)\in D_\epsilon(T)$. 
Set 
$$
\hat{v}_j(|x|,t):=(\partial_t^ju_*)(|x|,t)-F_N^j(|x|,t),
\qquad
\hat{u}_j(|x|,t):=U(|x|)\hat{v}_j(|x|,t).
$$ 
Since $F_N^j$ satisfies
$$
\frac{1}{\nu(r)r^{N-1}}\partial_r(\nu(r)r^{N-1}\partial_r F_N^j)=(\partial_t^{j+1} u_*)(r,t)
\quad\mbox{for $r>0$ and $t>0$},
$$
by Lemma~\ref{Lemma:2.2} and \eqref{eq:2.2} we have 
\begin{eqnarray}
\label{eq:3.11}
 & & \frac{1}{\nu(r)r^{N-1}}\partial_r(\nu(r)r^{N-1}\partial_r \hat{v}_j)=0\quad\mbox{for $r>0$ and $t>0$},\\
\label{eq:3.12}
 & & \limsup_{r\to 0} |\hat{v}_j(r,t)|<\infty\quad\mbox{for any $t>0$}.
\end{eqnarray}
It follows from \eqref{eq:3.11} that $\hat{u}_j$ satisfies (O) for any fixed $t>0$. 
On the other hand, 
since $U$ and $\tilde{U}$ are linearly independent solutions of (O), 
for any $t>0$, we can find constants $c_j(t)$ and $\hat{c}_j(t)$ such that 
$$
\hat{u}_j(r,t)=c_j(t)U(r)+\hat{c}_j(t)\tilde{U}(r)\quad\mbox{for $r>0$}. 
$$
This implies that 
$$
\hat{v}_j(r,t)=U(r)^{-1}\hat{u}_j(r,t)=c_j(t)+\hat{c}_j(t)U(r)^{-1}\tilde{U}(r)\quad\mbox{for $r>0$}. 
$$
Then, by \eqref{eq:1.2} and \eqref{eq:3.12} we have $\hat{c}_j(t)=0$ 
and see that $\hat{v}_j(r,t)\equiv c_j(t)$ for $r\ge 0$. 
Therefore we have
$$
(\partial_t^ju_*)(|x|,t)=c_j(t)+F_N^j(|x|,t)
\quad\mbox{in}\quad{\bf R}^N\times(0,\infty),
\qquad
c_j(t)=(\partial_t^ju_*)(0,t).
$$
Thus assertion~(ii) follows, and the proof of Proposition~\ref{Proposition:3.1} is complete. 
$\Box$
\section{Large time behavior of radially symmetric solutions}
In this section, under condition~(V), 
we study the large time behavior of radially symmetric solution~$u=u(|x|,t)$ of \eqref{eq:1.1} 
and prove Theorems~\ref{Theorem:1.1}--\ref{Theorem:1.3}.

Let $A$, $d$ and $w$ be as in Section~1.1.  
Set $U_d(r):=r^{-A}U(r)$ and $\nu_d:=U_d^2$. 
By \eqref{eq:1.2}, \eqref{eq:1.4} and Lemma~\ref{Lemma:2.4} we have:
\begin{eqnarray}
\notag
 & & \frac{1}{r^{d-1}}\partial_r\left(r^{d-1}\partial_rU_d\right)-V_{\lambda_2}(r)U_d=0\quad\mbox{in}\quad(0,\infty);\\
\notag
 & & U_d(r)\sim r^{A^+(\lambda_1)-A}=r^{A_d^+(\lambda)}\quad\mbox{as}\quad r\to 0;\\
\label{eq:4.1}
 & & U_d(r)\sim c_*,\quad U_d'(r)=O(r^{-1-\delta})\quad\mbox{as $r\to\infty$ in cases (S) and (C)}.\qquad\quad
\end{eqnarray}
Here $c_*$ is as in \eqref{eq:1.4}, $\lambda:=\lambda_1-\lambda_2$ and 
$$
A_d^+(\lambda):=\frac{-(d-2)+\sqrt{(d-2)^2+4\lambda}}{2}.
$$
Furthermore, similarly to Lemma~\ref{Lemma:2.4}, we see that 
\begin{equation}
\label{eq:4.2}
U_d(r)\sim c_*\log r,\quad U_d'(r)=O(r^{-1})\quad\mbox{as $r\to\infty$ in case of ($\mbox{S}_*$)}.
\end{equation}
Then the function $F_N^j$ given in Proposition~\ref{Proposition:3.1} satisfies  
\begin{equation}
\label{eq:4.3}
F_N^j(r,t)=
F_d^j(r,t):=\int_0^{r} s^{1-d}[\nu_d(s)]^{-1}\left(\int_0^s \tau^{d-1}\nu_d(\tau)(\partial_t^{j+1} u_*)(\tau,t)\,d\tau\right)\,ds,
\end{equation}
where $j\in\{0,1,2,\dots\}$. 
Furthermore,  
it follows from \eqref{eq:3.7} that 
\begin{equation}
\label{eq:4.4}
\int_0^s \tau^{d-1}\nu_d(\tau)\,d\tau\le Cs^d\nu_d(s),\quad s\ge 0.
\end{equation} 

Assume the same conditions as in Theorem~\ref{Theorem:1.1}. 
Let $\theta$ be the constant given in condition~(V) and set
$$
\theta_*=\frac{\theta}{4(2+\theta)}\in\left(0,\frac{\theta}{8}\right). 
$$
Since $\tilde{V}(\xi,s)=e^s V_{\lambda_2}(e^{\frac{s}{2}}\xi)$, 
it follows from (V)~(ii) that  
\begin{equation}
\label{eq:4.5}
|\tilde{V}(\xi,s)|\le C\xi^{-2}|e^{\frac{s}{2}}\xi|^{-\theta}
\le C\exp\left[-\frac{\theta}{2}s+(2+\theta)\theta_*s\right]=Ce^{-\frac{\theta}{4}s}
\end{equation}
for $\xi\in(e^{-\theta_*s},\infty)$ and $s>0$. 
Let $\delta$ be as in Lemma~\ref{Lemma:2.4}. 
Then, taking a sufficiently small $\theta>0$ if necessary, 
we have
\begin{equation}
\label{eq:4.6}
0<\theta<\min\{1,d,d^{-1}\},\qquad
\sigma:=\left(\frac{1}{2}-\theta_*\right)(1+\delta)-\frac{1}{2}>\theta_*>0. 
\end{equation}

We prepare some lemmas on estimates of $w$.
\begin{lemma}
\label{Lemma:4.1}
Let $\|\varphi_*\|_{L^2({\bf R}^N,\,\nu e^{|x|^2/4}\,dx)}=1$. 
Assume the same conditions as in Theorem~{\rm\ref{Theorem:1.1}}. 
Then 
\begin{itemize}
  \item[{\rm (i)}] 
  $\displaystyle{\sup_{s>0}\,e^{-\frac{d}{4}s}\|w(s)\|_{L^2({\bf R}_+,\rho_d\,d\xi)}<\infty}$:
  \item[{\rm (ii)}]  
  Assume that 
  \begin{equation}
  \label{eq:4.7}
  \sup_{s>0}\,e^{\gamma s}\|w(s)\|_{L^2({\bf R}_+,\rho_d\,d\xi)}<\infty
  \end{equation}
  for some $\gamma\ge -d/4$. 
  Then
  \begin{eqnarray}
  \label{eq:4.8}
   & & w(e^{-\theta_*s},s)=O(e^{-\gamma s}),\\
  \label{eq:4.9}
   & & (\partial_\xi w)(e^{-\theta_*s},s)=O(e^{-\gamma s-\theta_*s}),\\
  \label{eq:4.10}
   & & \int_0^{e^{-\theta_*s}}|w(\xi,s)|^2\rho_d\,d\xi=O(e^{-2\gamma s-d\theta_*s}),
  \end{eqnarray}
  for all sufficiently large $s>0$. 
\end{itemize}
\end{lemma}
{\bf Proof.}
Since 
\begin{equation}
\label{eq:4.11}
\begin{split}
w(\xi,s)=(1+t)^{\frac{d}{2}}r^{-A}u(r,t)
=(1+t)^{\frac{d}{2}}r^{-A}U(r)u_*(r,t)=(1+t)^{\frac{d}{2}}U_d(r)u_*(r,t) & \\
\mbox{with $\xi=(1+t)^{-\frac{1}{2}}r$ and $s=\log(1+t)$}, & 
\end{split}
\end{equation}
it follows from Lemma~\ref{Lemma:2.1}~(ii) that 
\begin{equation}
\label{eq:4.12}
\begin{split}
\|w(s)\|_{L^2({\bf R}_+,\rho_d\,d\xi)}^2
 & =(1+t)^{\frac{d}{2}}
 \int_0^\infty |u_*(r,t)|^2U(r)^2r^{N-1}e^{\frac{r^2}{4(1+t)}}\,dr\\
 & =(1+t)^{\frac{d}{2}}|{\bf S}^{N-1}|^{-1}\|u_*(t)\|_{L^2({\bf R}^N,e^{|x|^2/4(1+t)}\nu\,dx)}^2\\
 & \le (1+t)^{\frac{d}{2}}|{\bf S}^{N-1}|^{-1}\|\varphi_*\|_{L^2({\bf R}^N,e^{|x|^2/4}\nu\,dx)}^2\\
 & =(1+t)^{\frac{d}{2}}|{\bf S}^{N-1}|^{-1}\|\varphi\|_{L^2({\bf R}_+,e^{|x|^2/4}\,dx)}^2<\infty
\end{split}
\end{equation}
for $s>0$ and $t>0$ with $s=\log(1+t)$, 
where $|{\bf S}^{N-1}|$ is the volume of $(N-1)$-dimensional unit sphere, that is 
$|{\bf S}^{N-1}|=2\pi^{\frac{N}{2}}/\Gamma(N/2)$. 
Thus assertion~(i) follows. 

We prove assertion~(ii). 
It follows from \eqref{eq:4.12} that 
\begin{equation}
\label{eq:4.13}
\|w(s)\|_{L^2({\bf R}_+,\rho_d\,d\xi)}
\ge (1+t)^{\frac{d}{4}}|{\bf S}^{N-1}|^{-1/2}\|u_*(t)\|_{L^2({\bf R}^N,\,\nu\,dx)}
\end{equation}
for $s>0$ and $t>0$ with $s=\log(1+t)$. 
Assume \eqref{eq:4.7} for some $\gamma\ge -d/4$. 
Then 
$$
\sup_{t>0}\,(1+t)^{\gamma+\frac{d}{4}}\|u_*(t)\|_{L^2({\bf R}^N,\,\nu\,dx)}<\infty.
$$
Applying Proposition~\ref{Proposition:3.1} with $D=\gamma+d/4$ and $D'=0$, 
we obtain
\begin{equation}
\label{eq:4.14}
u_*(|x|,t)=u_*(0,t)+F_N^0(|x|,t)\quad\mbox{in}\quad{\bf R}^N\times(0,\infty). 
\end{equation}
Furthermore, for any $T>0$ and any sufficiently small $\epsilon>0$, 
\begin{equation}
\label{eq:4.15}
\begin{split}
 & |u_*(|x|,t)|\le Ct^{-\gamma-\frac{d}{2}},\\
 & |F_N^0(|x|,t)|\le Ct^{-\gamma-\frac{d}{2}-1}|x|^2\le C\epsilon^2t^{-\gamma-\frac{d}{2}},\\
 & |(\partial_r F_N^0)(|x|,t)|\le Ct^{-\gamma-\frac{d}{2}-1}|x|\le C\epsilon t^{-\gamma-\frac{d}{2}-\frac{1}{2}},
\end{split}
\end{equation}
for $(x,t)\in D_\epsilon(T)$.  
Then, by \eqref{eq:4.1}, \eqref{eq:4.11} and \eqref{eq:4.15} we have
$w(e^{-\theta_*s},s)=O(e^{-\gamma s})$
for all sufficiently large $s>0$. 
Furthermore, 
\begin{equation*}
\begin{split}
(\partial_\xi w)(e^{-\theta_*s},s)
 & =(1+t)^{\frac{d+1}{2}}[U_d'(r)u_*(r,t)+U_d(r)(\partial_r u_*)(r,t)]\\
 & =(1+t)^{\frac{d+1}{2}}[O(r^{-1-\delta})u_*(r,t)+(c_*+o(1))(\partial_r F_N^0)(r,t)]\\
 & = (1+t)^{-\gamma+\frac{1}{2}}O(r^{-1-\delta})+(1+t)^{-\gamma-\frac{1}{2}}O(r)\\
 & =O(e^{-\gamma s-\sigma s})+O(e^{-\gamma s-\theta_*s})=O(e^{-\gamma s-\theta_*s})
\end{split}
\end{equation*}
for all sufficiently large $s>0$, 
where $r=e^{\frac{1}{2}s-\theta_*s}$, $s=\log(1+t)$ and $\sigma$ is as in \eqref{eq:4.6}. 
So we have \eqref{eq:4.8} and \eqref{eq:4.9}. 

On the other hand, 
by \eqref{eq:4.1}, \eqref{eq:4.4}, \eqref{eq:4.11}, \eqref{eq:4.14} and \eqref{eq:4.15} we have
\begin{equation*}
\begin{split}
\int_0^{e^{-\theta_*s}}|w(\xi,s)|^2\rho_d\,d\xi
 & =(1+t)^{\frac{d}{2}}\int_0^{(1+t)^{\frac{1}{2}-\theta_*}} |u_*(r,t)|^2U_d(r)^2r^{d-1}e^{\frac{r^2}{4(1+t)}}\,dr\\
 & \le Ct^{-2\gamma-\frac{d}{2}}\int_0^{(1+t)^{\frac{1}{2}-\theta_*}} \nu_d(r)r^{d-1}\,dr\\
  & \le Ct^{-2\gamma-d\theta_*}U_d((1+t)^{\frac{1}{2}-\theta_*})^2=O(e^{-2\gamma s-d\theta_*s})
\end{split}
\end{equation*}
for all sufficiently large $s>0$ and $t>0$ with $s=\log(1+t)$.  
This implies \eqref{eq:4.10}. 
Thus assertion~(ii) follows, and the proof is complete. 
$\Box$
\begin{lemma}
\label{Lemma:4.2}
Assume the same conditions as in Lemma~{\rm\ref{Lemma:4.1}}. 
Then
\begin{equation}
\label{eq:4.16}
\sup_{s>0}\,\|w(s)\|_{L^2({\bf R}_+,\rho_d\,d\xi)}<\infty.
\end{equation}
\end{lemma}
{\bf Proof.}
Assume that \eqref{eq:4.7} holds for some $\gamma\ge -d/4$. 
Let $I(s):=(e^{-\theta_*s},\infty)$. 
It follows from \eqref{eq:1.16} that 
\begin{equation*}
\begin{split}
 & \frac{d}{ds}\int_{I(s)}|w(\xi,s)|^2\rho_d\,d\xi
 =2\int_{I(s)} w(\partial_s w)\rho_d\,d\xi+\theta_*e^{-\theta_*s}|w(e^{-\theta_*s},s)|^2\rho_d(e^{-\theta_*s})\\
 & =2\int_{I(s)} w\partial_\xi(\rho_d\partial_\xi w)\,d\xi+d\int_{I(s)}|w(\xi,s)|^2\rho_d\,d\xi\\
 & \qquad\qquad
 -2\int_{I(s)} \tilde{V}w^2\rho_d\,d\xi+\theta_*e^{-\theta_*s}|w(e^{-\theta_*s},s)|^2\rho_d(e^{-\theta_*s})\\
 & =-2w(e^{-\theta_*s},s)\rho_d(e^{-\theta_*s})(\partial_\xi w)(e^{-\theta_*s},s)
 -2\int_{I(s)} |\partial_\xi w|^2\rho_d\,d\xi+d\int_{I(s)} |w(\xi,s)|^2\rho_d\,d\xi\\
  & \qquad\qquad
  -2\int_{I(s)} \tilde{V}w^2\rho_d\,d\xi+\theta_*e^{-\theta_*s}|w(e^{-\theta_*s},s)|^2\rho_d(e^{-\theta_*s})
\end{split}
\end{equation*}
for $s>0$. 
This together with Lemma~\ref{Lemma:4.1} and \eqref{eq:4.5} implies that 
\begin{equation}
\label{eq:4.17}
\begin{split}
\frac{d}{ds}\int_{I(s)}
|w(\xi,s)|^2\rho_d\,d\xi
 & \le-2\int_{I(s)}|\partial_\xi w|^2\rho_d\,d\xi+d\int_{I(s)} |w(\xi,s)|^2\rho_d\,d\xi\\
 & \qquad
+Ce^{-\frac{\theta}{4}s}\int_{I(s)} |w(\xi,s)|^2\rho_d\,d\xi+O(e^{-2\gamma s}e^{-d\theta_*s})
\end{split}
\end{equation}
for all sufficiently large $s>0$.  

Set
\begin{equation}
\label{eq:4.18}
\hat{w}(\xi,s):=
\left\{
\begin{array}{ll}
w(\xi,s) & \mbox{if}\quad\xi\ge e^{-\theta_*s},\vspace{3pt}\\
w(e^{-\theta_*s},s) & \mbox{if}\quad 0\le\xi<e^{-\theta_*s}.
\end{array}
\right.
\end{equation}
It follows from Lemmas~\ref{Lemma:2.5} and \ref{Lemma:4.1} that
\begin{equation}
\label{eq:4.19}
\begin{split}
 & -2\int_{I(s)} |(\partial_\xi w)(\xi,s)|^2\rho_d\,d\xi+d\int_{I(s)}|w(\xi,s)|^2\rho_d\,d\xi\\
 & =-2\int_0^\infty|(\partial_\xi \hat{w})(\xi,s)|^2\rho_d\,d\xi+d\int_0^\infty |\hat{w}(\xi,s)|^2\rho_d\,d\xi
+O(e^{-d\theta_*s}w(e^{-\theta_*s},s)^2)\\
 & \le -2\mu_0\int_0^\infty |\hat{w}(\xi,s)|^2\rho_d\,d\xi+O(e^{-d\theta_*s}w(e^{-\theta_*s},s)^2)
 =O(e^{-2\gamma s-d\theta_*s})
\end{split}
\end{equation}
for all sufficiently large $s>0$. 
This together with \eqref{eq:4.17} implies that 
\begin{equation}
\label{eq:4.20}
\frac{d}{ds}\int_{I(s)} |w(\xi,s)|^2\rho_d\,d\xi
\le Ce^{-\frac{\theta}{4}s}\int_{I(s)} |w(\xi,s)|^2\rho_d\,d\xi+O(e^{-2\gamma s}e^{-d\theta_*s})
\end{equation}
for all sufficiently large $s>0$. 

On the other hand, 
by Lemma~\ref{Lemma:4.1}~(i) 
we see that \eqref{eq:4.7} holds with $\gamma=-d/4$. 
Without loss of generality, we can find $j\in\{0,1,2,\dots\}$ such that 
\begin{equation}
\label{eq:4.21}
j\theta_*<\frac{1}{2}<(j+1)\theta_*. 
\end{equation}
Since $\theta_*<1/4$, 
applying \eqref{eq:4.20} with $\gamma=-d/4$, we have 
$$
\int_{I(s)} |w(\xi,s)|^2\rho_d\,d\xi=O(e^{\frac{d}{2}s-d\theta_*s})
$$
for all sufficiently large $s>0$. 
This together with Lemma~\ref{Lemma:4.1} implies that 
\begin{equation}
\label{eq:4.22}
\sup_{s>0}\,e^{-\frac{d}{4}s+\frac{d}{2}\theta_*s}\|w(s)\|_{L^2({\bf R}_+,\rho_d\,d\xi)}<\infty. 
\end{equation}
By \eqref{eq:4.21}, 
repeating this arguments, we obtain 
$$
\sup_{s>0}\,e^{-\frac{d}{4}s+\frac{jd}{2}\theta_*s}\|w(s)\|_{L^2({\bf R}_+,\rho_d\,d\xi)}<\infty. 
$$
Applying \eqref{eq:4.20} with $\gamma=-d/4-jd\theta_*/2$ again, 
by \eqref{eq:4.21} we have 
$$
\sup_{s\ge 1}\int_{I(s)} |w(\xi,s)|^2\rho_d\,d\xi<\infty. 
$$
Then, similarly to \eqref{eq:4.22}, 
we obtain \eqref{eq:4.16}. 
Thus Lemma~\ref{Lemma:4.2} follows. 
$\Box$
\vspace{5pt}

Combining Lemma~\ref{Lemma:4.2} with assertion~(ii) of Lemma~\ref{Lemma:4.1}, 
we have: 
\begin{lemma}
\label{Lemma:4.3}
Assume the same conditions as in Lemma~{\rm\ref{Lemma:4.1}}. 
Let $\hat{w}$ be as in \eqref{eq:4.18}. 
Then 
\begin{equation*}
\begin{split}
 & \sup_{s\ge 1}\,|w(e^{-\theta_*s},s)|<\infty,
 \qquad\quad
 \sup_{s\ge 1}\,e^{\theta_*s}|(\partial_\xi w)(e^{-\theta_*s},s)|<\infty,\\
 & \sup_{s\ge 1}\,\|\hat{w}(s)\|_{L^2({\bf R}_+,\rho_d\,d\xi)}<\infty,
 \quad
 \sup_{s\ge 1}\,
 e^{d\theta_*s}\int_0^{e^{-\theta_*s}}|w(\xi,s)|^2\rho_d\,d\xi<\infty.
\end{split}
\end{equation*}
\end{lemma}

Next we study the large time behavior of $\hat{w}$ and prove the following proposition. 
\begin{proposition}
\label{Proposition:4.1}
Let $\|\varphi_*\|_{L^2({\bf R}^N,\,\nu e^{|x|^2/4}\,dx)}=1$. 
Assume the same conditions as in Theorem~{\rm\ref{Theorem:1.1}}. 
Let $\hat{w}$ be as in \eqref{eq:4.18}. Set 
$$
a(s):=\int_0^\infty\hat{w}\psi_d\rho_d\,d\xi
=c_d\int_0^\infty\hat{w}(\xi,s)\xi^{d-1}\,d\xi.
$$
Then 
$\|\hat{w}-a(s)\psi_d\|_{L^2({\bf R}_+,\rho_d\,d\xi)}=O(e^{-\theta' s})$ 
as $s\to\infty$, where $\theta':=\min\{d\theta_*/2,\theta/8\}$. 
\end{proposition}
For the proof of Proposition~\ref{Proposition:4.1}, 
we prepare the following lemma.
\begin{lemma}
\label{Lemma:4.4}
Assume the same conditions as in Proposition~{\rm\ref{Proposition:4.1}}. 
Then 
\begin{eqnarray}
\label{eq:4.23}
 & & \sup_{s\ge 1}\,|a(s)|<\infty,\\
\label{eq:4.24}
 & & \sup_{s\ge 1}\,e^{2\theta's}|a'(s)|<\infty,\\
\label{eq:4.25}
 & & \sup_{s\ge 1}\left|\frac{d}{ds}w(e^{-\theta_*s},s)\right|<\infty.
\end{eqnarray}
\end{lemma}
{\bf Proof.}
It follows from Lemma~\ref{Lemma:4.3} that 
$$
\sup_{s\ge 1}\,|a(s)|
\le\sup_{s\ge 1}\,\|\hat{w}\|_{L^2({\bf R}_+,\rho_d\,d\xi)}\|\psi_d\|_{L^2({\bf R}_+,\rho_d\,d\xi)}<\infty. 
$$
So we have \eqref{eq:4.23}. 

We prove \eqref{eq:4.25}.
By Proposition~\ref{Proposition:3.1}~(ii) and \eqref{eq:4.11} 
we have
$$
w(e^{-\theta_*s},s)=e^{\frac{d}{2}s}U_d(r(s))u_*(r(s),t(s))
=e^{\frac{d}{2}s}U_d(r(s))[u_*(0,t(s))+F_N^0(r(s),t(s))]
$$
for $s>0$, where $r(s)=e^{\frac{1}{2}s-\theta_*s}$ and $t(s)=e^s-1$. 
Then
\begin{equation}
\label{eq:4.26}
\begin{split}
\frac{d}{ds}w(e^{-\theta_*s},s)
 & =\frac{d}{2}w(e^{-\theta_*s},s)+\frac{U_d'(r(s))}{U_d(r(s))}r'(s)w(e^{-\theta_*s},s)\\
 & \quad
 +e^{\frac{d}{2}s}U_d(r(s))(\partial_t u_*)(0,t(s))t'(s)\\
 & \quad
+e^{\frac{d}{2}s}U_d(r(s))[(\partial_r F_N^0)(r(s),t(s))r'(s)+(\partial_tF_N^0)(r(s),t(s))t'(s)]
\end{split}
\end{equation}
for $s>0$. 
It follows from \eqref{eq:4.1} that 
\begin{equation}
\label{eq:4.27}
U_d(r(s))\sim c_*,
\qquad
U_d'(r(s))r'(s)=O(r(s)^{-1-\delta}r'(s))=O(e^{-\delta\left(\frac{1}{2}-\theta_*\right)s}),
\end{equation}
for all sufficiently large $s>0$. 
On the other hand, 
by Lemma~\ref{Lemma:4.2} and \eqref{eq:4.13} we have 
\begin{equation}
\label{eq:4.28}
\sup_{t>0}\,(1+t)^{\frac{d}{4}}\,\|u_*(t)\|_{L^2({\bf R}^N,\nu\,dx)}<\infty. 
\end{equation}
Then we apply Proposition~\ref{Proposition:3.1} with $D=d/4$ and $D'=0$ 
to obtain
\begin{equation}
\label{eq:4.29}
\begin{split}
 & (\partial_t u_*)(0,t(s))=O(e^{-\frac{d}{2}s-s}),\\
 & (\partial_r F_N^0)(r(s),t(s))=O(e^{-\frac{d}{2}s-s}r(s)),\\
 & (\partial_tF_N^0)(r(s),t(s))=F_N^1(r(s),t(s))=O(e^{-\frac{d}{2}s-2s}r(s)^2),
\end{split}
\end{equation}
for all sufficiently large $s>0$. 
By Lemma~\ref{Lemma:4.3}, \eqref{eq:4.26}, \eqref{eq:4.27} and \eqref{eq:4.29} we have \eqref{eq:4.25}. 
Furthermore, 
by Lemma~\ref{Lemma:2.5}, Lemma~\ref{Lemma:4.3}, \eqref{eq:4.5}, \eqref{eq:4.16} and \eqref{eq:4.25} 
we obtain
\begin{equation*}
\begin{split}
 & a'(s)=\frac{c_d}{d}e^{-d\theta_*s}[w(e^{-\theta_*s},s)]'
 -c_d\theta_*e^{-d\theta_*s}w(e^{-\theta_*s},s)\\
 & \qquad\qquad\qquad\qquad
 +c_d\theta_*e^{-d\theta_*s}w(e^{-\theta_*s},s)
 +\int_{I(s)} \partial_sw\psi_d\rho_d\,d\xi\\
 & =O(e^{-d\theta_*s})+ 
 \int_{I(s)} \partial_\xi(\rho_d\partial_\xi w)\psi_d\,d\xi
+\frac{d}{2}\int_{I(s)} w\psi_d\rho_d\,d\xi-\int_{I(s)}\tilde{V}w\psi_d\rho_d\,d\xi\\
 & =O(e^{-d\theta_*s})+\int_{I(s)} w\partial_\xi(\rho_d\partial_\xi\psi_d)\,d\xi
 +\frac{d}{2}\int_{I(s)} w\psi_d\rho_d\,d\xi+O(e^{-\frac{\theta}{4}s})\\
 & =O(e^{-d\theta_*s})+O(e^{-\frac{\theta}{4}s})
 \end{split}
\end{equation*}
for all sufficiently large $s>0$. 
This implies \eqref{eq:4.24}. Thus Lemma~\ref{Lemma:4.4} follows.
$\Box$
\vspace{5pt}
\newline
{\bf Proof of Proposition~\ref{Proposition:4.1}.}
Set $\tilde{w}(\xi,s):=\hat{w}(\xi,s)-a(s)\psi_d(\xi)$. 
It follows from Lemma~\ref{Lemma:2.5} and \eqref{eq:1.16} that
$$
\partial_s \tilde{w}=\partial_s\hat{w}-a'(s)\psi_d
=-{\mathcal L}_d\hat{w}-\tilde{V}\hat{w}-a'(s)\psi_d
=-{\mathcal L}_d\tilde{w}-\tilde{V}\hat{w}-a'(s)\psi_d
$$
for $\xi\in I(s)$ and $s>0$. 
By Lemma~\ref{Lemma:4.3}, Lemma~\ref{Lemma:4.4} and \eqref{eq:4.5} we have
\begin{equation}
\label{eq:4.30}
\begin{split}
 & \frac{d}{ds}\int_{I(s)} |\tilde{w}(\xi,s)|^2\rho_d\,d\xi
 =2\int_{I(s)} \tilde{w}(\partial_s \tilde{w})\rho_d\,d\xi+\theta_*e^{-\theta_*s}|\tilde{w}(e^{-\theta_*s},s)|^2\rho_d(e^{-\theta_*s})\\
 & =2\int_{I(s)} \left[\tilde{w}\partial_\xi(\rho_d\partial_\xi \tilde{w})
 +\frac{d}{2}\tilde{w}^2\rho_d-\tilde{V}\hat{w}\tilde{w}\rho_d-a'(s)\psi_d\tilde{w}\rho_d\right]\,d\xi
 +O(e^{-d\theta_*s})\\
 & =-2\tilde{w}(e^{-\theta_*s},s)\rho_d(e^{-\theta_*s})(\partial_\xi \tilde{w})(e^{-\theta_*s},s)\\
 & \qquad\
 -2\int_{I(s)} |(\partial_\xi \tilde{w})(\xi,s)|^2\rho_d\,d\xi+d\int_{I(s)} |\tilde{w}(\xi,s)|^2\rho_d\,d\xi\\
 & \qquad\qquad
  -2\int_{I(s)}\tilde{V}\hat{w}\tilde{w}\rho_d\,d\xi
  -2a'(s)\int_{I(s)} \tilde{w}\psi_d\rho_d\,d\xi
  +O(e^{-d\theta_*s})\\
  & =-2\int_{I(s)}
  |(\partial_\xi \tilde{w})(\xi,s)|^2\rho_d\,d\xi+d\int_{I(s)} |\tilde{w}(\xi,s)|^2\rho_d\,d\xi+O(e^{-d\theta_*s})+O(e^{-\frac{\theta}{4}s})
\end{split}
\end{equation}
for all sufficiently large $s>0$. 
Furthermore, 
similarly to \eqref{eq:4.19}, 
by Lemmas~\ref{Lemma:2.5}, \ref{Lemma:4.3} and \ref{Lemma:4.4} 
we obtain 
\begin{equation}
\label{eq:4.31}
\begin{split}
 & \int_{I(s)} |(\partial_\xi \tilde{w})(\xi,s)|^2\rho_d\,d\xi-\frac{d}{2}\int_{I(s)} |\tilde{w}(\xi,s)|^2\rho_d\,d\xi\\
 & \ge\int_0^\infty|(\partial_\xi \tilde{w})(\xi,s)|^2\rho_d\,d\xi-\frac{d}{2}\int_0^\infty |\tilde{w}(\xi,s)|^2\rho_d\,d\xi
 -a(s)^2\int_0^{e^{-\theta_*s}}|\partial_\xi\psi_d|^2\rho_d\,d\xi\\
 & =\int_0^\infty|(\partial_\xi \tilde{w})(\xi,s)|^2\rho_d\,d\xi-\frac{d}{2}\int_0^\infty |\tilde{w}(\xi,s)|^2\rho_d\,d\xi+O(e^{-(d+2)\theta_*s})\\
 & \ge\int_0^\infty|\tilde{w}(\xi,s)|^2\rho_d\,d\xi+O(e^{-(d+2)\theta_*s})
 \ge\int_{I(s)} |\tilde{w}(\xi,s)|^2\rho_d\,d\xi+O(e^{-(d+2)\theta_*s})
\end{split}
\end{equation}
for all sufficiently large $s>0$. 
Therefore we deduce from \eqref{eq:4.30} and \eqref{eq:4.31} that
\begin{equation}
\label{eq:4.32}
\begin{split}
\frac{d}{ds}\int_{I(s)} |\tilde{w}(\xi,s)|^2\rho_d\,d\xi
 & \le -2\int_{I(s)} |\tilde{w}(\xi,s)|^2\rho_d\,d\xi+O(e^{-d\theta_*s})+O(e^{-\frac{\theta}{4}s})\\
 & =-2\int_{I(s)} |\tilde{w}(\xi,s)|^2\rho_d\,d\xi+O(e^{-2\theta' s})
\end{split}
\end{equation}
for all sufficiently large $s>0$. 
Since $d\theta_*<d\theta<1$ (see \eqref{eq:4.6}), 
by \eqref{eq:4.32} we have 
\begin{equation}
\label{eq:4.33}
\int_{I(s)} |\tilde{w}(\xi,s)|^2\rho_d\,d\xi=O(e^{-2\theta' s})
\end{equation}
for all sufficiently large $s>0$. 
Combining \eqref{eq:4.33} with Lemmas~\ref{Lemma:4.3} and \ref{Lemma:4.4}, 
we obtain 
$$
\int_0^\infty |\tilde{w}(\xi,s)|^2\rho_d\,d\xi=O(e^{-2\theta' s})
$$
for all sufficiently large $s>0$. 
Thus Proposition~\ref{Proposition:4.1} follows. 
$\Box$
\begin{proposition}
\label{Proposition:4.2}
Let $\|\varphi_*\|_{L^2({\bf R}^N,\,\nu e^{|x|^2/4}\,dx)}=1$. 
Assume the same conditions as in Theorem~{\rm\ref{Theorem:1.1}}. 
Then 
\begin{eqnarray}
\label{eq:4.34}
 & & |a(s)-m(\varphi)|=O(e^{-2\theta's}),\\
\label{eq:4.35}
 & & \|\hat{w}(s)-m(\varphi)\psi_d\|_{L^2({\bf R}_+,\rho_d\,d\xi)}=O(e^{-\theta's}),
\end{eqnarray}
for all sufficiently large $s>0$. 
Furthermore, if $m(\varphi)=0$, then 
\begin{equation}
\label{eq:4.36}
\|\hat{w}(s)\|_{L^2({\bf R}_+,\rho_d\,d\xi)}=O(e^{-s}),
\quad
\|w(s)\|_{L^2({\bf R}_+,\rho_d\,d\xi)}=O(e^{-s}),
\end{equation}
for all sufficiently large $s>0$. 
\end{proposition}
{\bf Proof.}
By \eqref{eq:4.24} we can find a constant $a_\infty$ such that 
\begin{equation}
\label{eq:4.37}
|a(s)-a_\infty|=O(e^{-2\theta's})\quad\mbox{as}\quad s\to\infty. 
\end{equation}
On the other hand, 
by Lemma~\ref{Lemma:4.3} we have 
\begin{equation}
\label{eq:4.38}
\left|\int_{I(s)^c}\hat{w}\xi^{d-1}\,d\xi\right|
\le\left(\int_{I(s)^c}\hat{w}^2\rho_d\,d\xi\right)^{1/2}\left(\int_{I(s)^c}\xi^{d-1}e^{-\frac{\xi^2}{4}}\,d\xi\right)^{1/2}
=O(e^{-d\theta_*s})
\end{equation}
for all sufficiently large $s>0$, 
where $I(s)^c:={\bf R}_+\setminus I(s)$. 
By Lemma~\ref{Lemma:4.3}, \eqref{eq:4.1}, \eqref{eq:4.11} and \eqref{eq:4.38}
we obtain 
\begin{equation}
\label{eq:4.39}
\begin{split}
a(s) & =c_d\int_{I(s)} w\xi^{d-1}\,d\xi+O(e^{-d\theta_*s})\\
 & =c_d\int_{(1+t)^{\frac{1}{2}-\theta_*}}^\infty u_*(r,t)U_d(r)r^{d-1}\,dr+O(e^{-d\theta_*s})\\
 & =\frac{c_d}{c_*}\int_{(1+t)^{\frac{1}{2}-\theta_*}}^\infty u_*(r,t) \nu_d(r)r^{d-1}\,dr+o(1)
\end{split}
\end{equation}
for all sufficiently large $s>0$ and $t>0$ with $s=\log(1+t)$. 

On the other hand, by \eqref{eq:4.28}
we apply Proposition~\ref{Proposition:3.1} with $D=d/4$ and $D'=0$ to obtain 
\begin{equation}
\label{eq:4.40}
\sup_{0\le r\le (1+t)^{\frac{1}{2}-\theta_*}}|u_*(r,t)|=O(t^{-\frac{d}{2}})
\end{equation}
for all sufficiently large $t>0$. 
Combining \eqref{eq:4.40} with \eqref{eq:4.4}, we see that 
\begin{equation}
\label{eq:4.41}
\begin{split}
\int_0^{(1+t)^{\frac{1}{2}-\theta_*}}u_*(r,t) \nu_d(r)r^{d-1}\,dr
 & =O(t^{-\frac{d}{2}})\int_0^{(1+t)^{\frac{1}{2}-\theta_*}}\nu_d(r)r^{d-1}\,dr\\
 & =O(t^{-\frac{d}{2}})O(t^{\frac{d}{2}-d\theta_*})
=O(t^{-d\theta_*})
\end{split}
\end{equation}
for all sufficiently large $t>0$.
Therefore, by \eqref{eq:4.39} and \eqref{eq:4.41} we obtain 
\begin{equation}
\label{eq:4.42}
a_\infty=\lim_{s\to\infty}a(s)=\lim_{t\to\infty}\frac{c_d}{c_*}\int_0^\infty u_*(r,t)\nu_d(r)r^{d-1}\,dr.
\end{equation}
On the other hand, since $u_*$ is a radial solution of problem~(P), 
we have 
\begin{equation}
\label{eq:4.43}
\int_0^\infty u_*(r,t)\nu_d(r)r^{d-1}\,dr
=\int_0^\infty u_*(r,t)\nu(r)r^{N-1}\,dr
=\int_0^\infty \varphi_*(r)\nu(r)r^{N-1}\,dr.
\end{equation}
We deduce from \eqref{eq:4.42} and \eqref{eq:4.43} that $a_\infty=m(\varphi)$. 
This together with \eqref{eq:4.37} implies \eqref{eq:4.34}. 
Furthermore, by Proposition~\ref{Proposition:4.1} and \eqref{eq:4.34} we have \eqref{eq:4.35}. 

It remains to prove \eqref{eq:4.36}. 
Assume that $m(\varphi)=0$. 
Then it follows from \eqref{eq:4.35} and Lemma~\ref{Lemma:4.3} that 
\begin{equation}
\label{eq:4.44}
\|\hat{w}(s)\|_{L^2({\bf R}_+,\rho_d\,d\xi)}=O(e^{-\theta's}),
\qquad
\|w(s)\|_{L^2({\bf R}_+,\rho_d\,d\xi)}=O(e^{-\theta's}),
\end{equation}
for all sufficiently large $s>0$. 
Applying the same argument as in the proof of \eqref{eq:4.32}, 
we see that
$$
\frac{d}{ds}\int_{I(s)} |\tilde{w}(\xi,s)|^2\rho_d\,d\xi
\le -\int_{I(s)} |\tilde{w}(\xi,s)|^2\rho_d\,d\xi+O(e^{-4\theta's})
$$
for all sufficiently large $s>0$. 
Furthermore, similarly to \eqref{eq:4.44}, 
we have 
$$
\|\hat{w}(s)\|_{L^2({\bf R}_+,\rho_d\,d\xi)}=O(e^{-2\theta's}),
\qquad
\|w(s)\|_{L^2({\bf R}_+,\rho_d\,d\xi)}=O(e^{-2\theta's}),
$$
for all sufficiently large $s>0$. 
Repeating this argument, we can find $\tilde{\theta}>1$ such that 
$$
\frac{d}{ds}\int_{I(s)} |\tilde{w}(\xi,s)|^2\rho_d\,d\xi
\le -\int_{I(s)} |\tilde{w}(\xi,s)|^2\rho_d\,d\xi+O(e^{-\tilde{\theta} s})
$$
for all sufficiently large $s>0$, instead of \eqref{eq:4.32}. 
This implies that 
$$
\|\hat{w}(s)\|_{L^2({\bf R}_+,\rho_d\,d\xi)}=O(e^{-s}),
\qquad
\|w(s)\|_{L^2({\bf R}_+,\rho_d\,d\xi)}=O(e^{-s}),
$$
for all sufficiently large $s>0$. 
Thus \eqref{eq:4.36} holds. 
Therefore the proof of Proposition~\ref{Proposition:4.2} is complete. 
$\Box$\vspace{5pt}

We are ready to complete the proof of Theorems~\ref{Theorem:1.1} and \ref{Theorem:1.2}.
\vspace{5pt}
\newline
{\bf Proof of Theorem~\ref{Theorem:1.1}.}
By the linearity of the operator $L$ 
it suffices to consider only the case 
\begin{equation}
\label{eq:4.45}
1=\|\varphi\|_{L^2({\bf R}^N,e^{|x|^2/4}\,dx)}=\|\varphi_*\|_{L^2({\bf R}^N,\,\nu e^{|x|^2/4}\,dx)}
=|{\bf S}^{N-1}|^{1/2}\|w(0)\|_{L^2({\bf R}_+,\,\rho_d\,d\xi)}.
\end{equation}
Let $R>1$. 
By Lemma~\ref{Lemma:4.2} 
we apply the parabolic regularity theorems (see e.g., \cite{LSU}) to \eqref{eq:1.16}. 
Then we can find $\alpha\in(0,1)$ such that
\begin{equation}
\label{eq:4.46}
\|w\|_{C^{2,\alpha;1,\alpha/2}(\{R^{-1}\le|y|\le R\}\times(S,\infty))}<\infty
\end{equation}
for any $R>1$ and $S>0$. 
Then, by Proposition~\ref{Proposition:4.2} and \eqref{eq:1.16} 
we apply the Ascoli-Arzel\`a theorem and the diagonal argument to obtain 
\begin{equation}
\label{eq:4.47}
\lim_{s\to\infty}
\left\|w(s)-m(\varphi)\psi_d\right\|_{C^2(\{R^{-1}\le|y|\le R\})}=0,
\quad
\lim_{s\to\infty}
\left\|(\partial_sw)(s)\right\|_{C^2(\{R^{-1}\le|y|\le R\})}=0.
\end{equation}
Furthermore, if $a_\infty=m(\varphi)=0$, then, 
similarly to \eqref{eq:4.46}, by \eqref{eq:4.36} we have
$$
\sup\,
\left\{\left|(\partial_\xi^\ell w)(\xi,s)\right|\,:\,
R^{-1}\le|y|\le R,\,s\ge S\right\}=O(e^{-s})
\quad\mbox{as}\quad s\to\infty
$$
for any $R>1$, where $\ell=0,1,2$.
These together with Proposition~\ref{Proposition:4.2} imply \eqref{eq:1.10} and \eqref{eq:1.12}. 
Thus Theorem~\ref{Theorem:1.1} follows.
$\Box$
\vspace{5pt}
\newline
{\bf Proof of Theorem~\ref{Theorem:1.2}.}
Similarly to the proof of Theorem~\ref{Theorem:1.1}, 
we can assume \eqref{eq:4.45} without loss of generality. 
Let $T>0$ and let $\epsilon$ be any sufficiently small positive constant. 
By Lemma~\ref{Lemma:4.2} and \eqref{eq:4.3}, 
applying Proposition~\ref{Proposition:3.1} with $D=d/4$ and $D'=0$, 
we obtain 
\begin{equation}
\label{eq:4.48}
(\partial_t^ju_*)(|x|,t)=(\partial_t^ju_*)(0,t)+F_d^j(|x|,t)\quad\mbox{in}\quad{\bf R}^N\times(0,\infty),
\end{equation}
where $j\in\{0,1,2,\dots\}$. Furthermore, 
\begin{equation}
\label{eq:4.49}
|F_d^j(r,t)|\le Ct^{-\frac{d}{2}-j-1}r^2,
\qquad
|(\partial_rF_d^j)(r,t)|\le Ct^{-\frac{d}{2}-j-1}r,
\end{equation}
for $0\le r\le\epsilon(1+t)^{1/2}$ and $t\ge T$.
Then it follows from \eqref{eq:4.48} and \eqref{eq:4.49} that 
\begin{equation}
\label{eq:4.50}
|(\partial_r u_*)(r,t)|\le C_2t^{-\frac{d}{2}-1}r
\end{equation}
for $0\le r\le\epsilon(1+t)^{1/2}$ and $t\ge T$.
Furthermore, 
by \eqref{eq:4.3} and \eqref{eq:4.48} we have 
\begin{equation}
\label{eq:4.51}
\begin{split}
F_d^0(r,t)
 & =\int_0^r s^{1-d}[\nu_d(s)]^{-1}\left(\int_0^s \tau^{d-1}\nu_d(\tau)(\partial_t u_*)(\tau,t)\,d\tau\right)\,ds\\
 & =\int_0^r s^{1-d}[\nu_d(s)]^{-1}\left(\int_0^s \tau^{d-1}\nu_d(\tau)
 \biggr[(\partial_tu_*)(0,t)+F_d^1(\tau,t)\biggr]\,d\tau\right)\,ds\\
 & =(\partial_tu_*)(0,t)F_d(r)+G_d(r,t)
\end{split}
\end{equation}
for $r\ge 0$ and $t>0$, 
where $F_d$ is given in Theorem~\ref{Theorem:1.2} and 
\begin{equation}
\label{eq:4.52}
G_d(r,t)=\int_0^r s^{1-d}[\nu_d(s)]^{-1}\left(\int_0^s \tau^{d-1}\nu_d(\tau)
F_d^1(\tau,t)\,d\tau\right)\,ds.
\end{equation}
Then \eqref{eq:1.13} holds. 
In addition, 
by \eqref{eq:4.4}, \eqref{eq:4.49} and \eqref{eq:4.52} we have 
\begin{equation}
\label{eq:4.53}
\begin{split}
|G_d(r,t)| & \le Ct^{-\frac{d}{2}-2}
\int_0^r s^{1-d}[\nu_d(s)]^{-1}\left(\int_0^s \tau^{d+1}\nu_d(\tau)\,d\tau\right)\,ds\\
 & \le Ct^{-\frac{d}{2}-2}
\int_0^r s^{1-d}[\nu_d(s)]^{-1}\cdot s^{d+2}\nu_d(s)\,ds
\le Ct^{-\frac{d}{2}-2}r^4
\end{split}
\end{equation}
for $0\le r\le\epsilon(1+t)^{1/2}$ and $t\ge T$. 
A similar argument with \eqref{eq:4.1} implies that 
$$
|(\partial_r^\ell G_d)(r,t)|\le Ct^{-\frac{d}{2}-2}r^{4-\ell}
$$
for $0\le r\le\epsilon(1+t)^{1/2}$ and $t\ge T$, where $\ell\in\{1,2\}$. 
Thus \eqref{eq:1.14} holds for $\ell\in\{0,1,2\}$. 

It remains to prove assertion~(b). 
By \eqref{eq:4.11} and \eqref{eq:4.48} 
we have 
\begin{equation}
\label{eq:4.54}
w(\xi,s)=(1+t)^{\frac{d}{2}}U_d(r)u_*(r,t)=(1+t)^{\frac{d}{2}}U_d(r)
 \left[u_*(0,t)+F_d^0(r,t)\right]
\end{equation}
for $\xi\in(0,\infty)$ and $s>0$ with $\xi=(1+t)^{-1/2}r$ and $s=\log(1+t)$. 
Let $0<\xi<\epsilon$. 
By \eqref{eq:4.1}, \eqref{eq:4.49} and \eqref{eq:4.54} we obtain 
\begin{equation}
\label{eq:4.55}
\left|w(\xi,s)
-(1+t)^{\frac{d}{2}}(c_*+o(1))u_*(0,t)\right|\le C\xi^2
\end{equation}
for all sufficiently large $s>0$ and $t>0$ with $s=\log(1+t)$ and $0<\xi<\epsilon$. 
On the other hand, it follows from \eqref{eq:4.47} that
\begin{equation}
\label{eq:4.56}
\lim_{s\to\infty}w(\xi,s)=c_dm(\varphi)e^{-\frac{\xi^2}{4}}.
\end{equation}
Then we deduce from \eqref{eq:4.55} and \eqref{eq:4.56} that 
\begin{equation}
\label{eq:4.57}
\lim_{t\to\infty}t^{\frac{d}{2}}u_*(0,t)=\frac{c_d}{c_*}m(\varphi).
\end{equation}
Furthermore, it follows from \eqref{eq:4.11} that
\begin{equation*}
\begin{split}
(\partial_s w)(\xi,s) & =\frac{d}{2}w(\xi,s)
+e^{\frac{(d+1)s}{2}}U_d'(e^{\frac{s}{2}}\xi)\frac{\xi}{2}u_*(e^{\frac{s}{2}}\xi,t)\\
 & +e^{\frac{(d+1)s}{2}}U_d(e^{\frac{s}{2}}\xi)\frac{\xi}{2}(\partial_ru_*)(e^{\frac{s}{2}}\xi,t)
+e^{\frac{ds}{2}+s}U_d(e^{\frac{s}{2}}\xi)(\partial_tu_*)(e^{\frac{s}{2}}\xi,e^s-1).
\end{split}
\end{equation*}
This together with \eqref{eq:4.1}, \eqref{eq:4.50} and \eqref{eq:4.56} implies that 
\begin{equation}
\label{eq:4.58}
\begin{split}
(\partial_s w)(\xi,s)
& =\frac{d}{2}w(\xi,s)+e^{\frac{s}{2}}\frac{U_d'(e^{\frac{s}{2}}\xi)}{U_d(e^{\frac{s}{2}}\xi)}\frac{\xi}{2}w(\xi,s)\\
& \qquad
+O(\xi^2)+e^{\frac{ds}{2}+s}(c_*+o(1))(\partial_tu_*)(e^{\frac{s}{2}}\xi,e^s-1)\\
& =\frac{d}{2}m(\varphi)c_de^{-\frac{\xi^2}{4}}+o(1)+O((e^{\frac{s}{2}}\xi)^{-\delta})\\
 & \qquad
+O(\xi^2)+e^{\frac{ds}{2}+s}(c_*+o(1))(\partial_tu_*)(e^{\frac{s}{2}}\xi,e^s-1)
\end{split}
\end{equation}
for all sufficiently large $s>0$.
On the other hand, by \eqref{eq:4.48} and \eqref{eq:4.49} we have 
\begin{equation}
\label{eq:4.59}
e^{\frac{ds}{2}+s}(\partial_tu_*)(e^{\frac{s}{2}}\xi,e^s-1)
=e^{\frac{ds}{2}+s}(\partial_t u_*)(0,e^s-1)+O(\xi^2)
\end{equation}
for all sufficiently large $s>0$. 
Therefore, by \eqref{eq:4.47}, \eqref{eq:4.58} and \eqref{eq:4.59} 
we obtain 
$$
\limsup_{s\to\infty}\left|(c_*+o(1))e^{\frac{ds}{2}+s}(\partial_t u_*)(0,e^s-1)+\frac{d}{2}c_dm(\varphi)\right|
\le C\xi^2.
$$
Since $0<\xi<\epsilon$, we deduce that 
$$
\limsup_{s\to\infty}\left|e^{\frac{ds}{2}+s}(\partial_t u_*)(0,e^s-1)+\frac{dc_d}{2c_*}m(\varphi)\right|=0. 
$$
This together with \eqref{eq:4.57} implies assertion~(b). 
Thus Theorem~\ref{Theorem:1.2} follows.
$\Box$
\vspace{5pt}

\noindent
{\bf Proof of Theorem~\ref{Theorem:1.3}.}
Similarly to the proof of Theorem~\ref{Theorem:1.1}, 
we can assume \eqref{eq:4.45} without loss of generality. 
Assume the same conditions as in Theorem~\ref{Theorem:1.3}. 
Then $d=2$ and $V_{\lambda_2}$ satisfies condition~(V) 
with $\lambda_1$ and $\lambda_2$ replaced by $\lambda_1-\lambda_2\,(\ge 0)$ and $0$, respectively. 
Applying a similar argument as in the proof of argument as in \cite[Proposition~3.1]{IK04}, 
we have 
\begin{equation}
\label{eq:4.60}
\lim_{s\to\infty}sw(\xi,s)
=\frac{1}{c_*}\left[\int_0^\infty w(r,0)U_d(r)r\,dr\right]e^{-\frac{\xi^2}{4}}
=2m(\varphi)\psi_d(\xi)
\end{equation}
in $L^2({\bf R}_+,\rho_2\,d\xi)\cap C^2(K)$, 
for any compact set $K$ in ${\bf R}^2\setminus\{0\}$.  
Furthermore, 
$$
\lim_{t\to\infty}t(\log t)^2u_*(0,t)=2\sqrt{2}c_*^{-1}m(\varphi),
\quad
\lim_{t\to\infty}t^2(\log t)^2(\partial_tu_*)(0,t)=-2\sqrt{2}c_*^{-1}m(\varphi). 
$$
On the other hand, similarly to \eqref{eq:4.48}, 
we have 
$$
(\partial_t^ju_*)(|x|,t)=(\partial_t^ju_*)(0,t)+F_2^j(|x|,t)\quad\mbox{in}\quad{\bf R}^N\times(0,\infty),
$$
where $j\in\{0,1,2,\dots\}$. 
It follows from \eqref{eq:4.60} that 
$$
\sup_{t>0}\,(1+t)^{\frac{d}{4}}\log(2+t)\,\|u_*(t)\|_{L^2({\bf R}^N,\nu e^{|x|^2/4(1+t)}\,dx)}<\infty.
$$
Let $T>0$ and $\epsilon$ be a sufficiently small positive constant. 
Then, by \eqref{eq:4.3} we apply Proposition~\ref{Proposition:3.1} with $D=d/4$ and $D'=1$ 
to obtain 
$$
\left|F_2^j(r,t)\right|\le Ct^{-1-j-1}[\log(2+t)]^{-2}r^2
$$
for $0\le r\le\epsilon(1+t)^{1/2}$ and $t\ge T$. 
Similarly to \eqref{eq:4.51} and \eqref{eq:4.52}, 
we have 
\begin{equation*}
\begin{split}
 & F_2^0(r,t)=(\partial_tu_*)(0,t)F_2(r)+G_2(r,t),\\
 & G_2(r,t)=\int_0^r s^{-1}[\nu_2(s)]^{-1}\left(\int_0^s \tau\nu_2(\tau)
F_2^1(\tau,t)\,d\tau\right)\,ds,
\end{split}
\end{equation*}
for $r\ge 0$ and $t>0$. 
Furthermore, similarly to \eqref{eq:4.53}, we obtain  
\begin{equation*}
\begin{split}
|G_2(r,t)| & \le Ct^{-3}[\log(2+t)]^{-2}
\int_0^r s^{-1}[\nu_2(s)]^{-1}\left(\int_0^s \tau^3\nu_2(\tau)\,d\tau\right)\,ds\\
 & \le Ct^{-3}[\log(2+t)]^{-2}r^4
\end{split}
\end{equation*}
for $0\le r\le\epsilon(1+t)^{1/2}$ and $t\ge T$. 
A similar argument with \eqref{eq:4.2} implies that 
$$
|(\partial_r^\ell G_2)(r,t)|\le Ct^{-3}[\log(2+t)]^{-2}r^{4-\ell},\quad \ell=1,2, 
$$
for $0\le r\le\epsilon(1+t)^{1/2}$ and $t\ge T$. 
So we see that \eqref{eq:1.15} holds for $\ell\in\{0,1,2\}$. 
Thus Theorem~\ref{Theorem:1.3} follows.
$\Box$
\section{Proof of Theorem~\ref{Theorem:1.4}}
We use the same notation as in Section~1.2. 
Let $m\in\{1,2,\dots\}$. 
Then 
$$
L_m:=-\Delta+V(|x|)+\frac{\omega_m}{|x|^2}
$$ 
is subcritical and 
problem~(O) corresponding to $L_m$ possesses a positive solution $U_m$ satisfying 
\begin{equation}
\label{eq:5.1}
U_m(r)\thicksim r^{A^+(\lambda_1+\omega_m)}\quad\mbox{as}\quad r\to+0,\quad
U_m(r)\thicksim c_m\,r^{A^+(\lambda_2+\omega_m)}\quad\mbox{as}\quad r\to\infty,
\end{equation}
for some positive constant $c_m$. 
Set
$$
u(x,t):=e^{-tL}\varphi,\qquad
u_m(x,t):=u(x,t)-\sum_{k=0}^{m-1}\sum_{i=1}^{\ell_k}
e^{-tL}\varphi^{k,i}.
$$
\begin{lemma}
\label{Lemma:5.1}
Let $m\in\{1,2,\dots\}$. Then there exists $C_1>0$ such that 
\begin{equation}
\label{eq:5.2}
\|u_m(t)\|_{L^2({\bf R}^N,e^{|x|^2/4(1+t)}\,dx)}
\le C_1t^{-\frac{d_m}{4}}\|u_m(0)\|_{L^2({\bf R}^N,e^{|x|^2/4}\,dx)}
\end{equation}
for $t>0$, where $d_m:=N+2A^+(\lambda_2+\omega_m)$. 
Furthermore, there exists $C_2>0$ such that 
\begin{equation}
\label{eq:5.3}
\left|\frac{u_m(x,t)}{U(\min\{|x|,\sqrt{t}\})}\right|\le C_2t^{-\frac{N+d_m}{4}}U(\sqrt{t})^{-1}\|u_m(0)\|_{L^2({\bf R}^N,e^{|x|^2/4}\,dx)}
\end{equation}
for $x\in{\bf R}^N$ and $t>0$. 
\end{lemma}
{\bf Proof.}
Let $m\in\{1,2,\dots\}$. 
The comparison principle implies that
$$
\left|\left[e^{-tL_k}\phi^{k,i}\right](x)\right|\le \left[e^{-tL_k}|\phi^{k,i}|\right](x)
\le \left[e^{-tL_m}|\phi^{k,i}|\right](x)
\quad\mbox{in}\quad{\bf R}^N\times(0,\infty)
$$
for $k\in\{m,m+1,\dots\}$ and $i\in\{1,\dots,\ell_k\}$. 
On the other hand, by Theorem~\ref{Theorem:1.1} and \eqref{eq:5.1} (see also \eqref{eq:4.28}) we have
$$
\left\|e^{-tL_m}|\phi^{k,i}|\right\|_{L^2({\bf R}^N,e^{|x|^2/4(1+t)}\,dx)}\le C(1+t)^{-\frac{d_m}{4}}\|\phi^{k,i}\|_{L^2({\bf R}^N,e^{|x|^2/4}\,dx)},
\quad t>0,
$$
for $k\in\{m,m+1,\dots\}$ and $i\in\{1,\dots,\ell_k\}$. 
These together with \eqref{eq:1.17} implies that 
\begin{equation*}
\begin{split}
 & \left\|e^{-tL}\varphi^{k,i}\right\|_{L^2({\bf R}^N,e^{|x|^2/4(1+t)}\,dx)}
=|{\bf S}^{N-1}|^{-\frac{1}{2}}\left\|e^{-tL_k}\phi^{k,i}\right\|_{L^2({\bf R}^N,e^{|x|^2/4(1+t)}\,dx)}\\
 & \qquad
\le C(1+t)^{-\frac{d_m}{4}}\|\phi^{k,i}\|_{L^2({\bf R}^N,e^{|x|^2/4}\,dx)}
\le C(1+t)^{-\frac{d_m}{4}}\|\varphi^{k,i}\|_{L^2({\bf R}^N,e^{|x|^2/4}\,dx)}
\end{split}
\end{equation*}
for $t>0$. 
Therefore we deduce from the orthogonality of $\{Q_{k,i}\}$ that
\begin{equation*}
\begin{split}
 & \|u_m(t)\|_{L^2({\bf R}^N,e^{|x|^2/4(1+t)}\,dx)}^2
=\sum_{k=m}^\infty\sum_{i=1}^{\ell_k}\left\|e^{-tL}\varphi^{k,i}\right\|_{L^2({\bf R}^N,e^{|x|^2/4(1+t)}\,dx)}^2\\
 & 
\le C(1+t)^{-\frac{d_m}{2}}\sum_{k=m}^\infty\sum_{i=1}^{\ell_k}\|\varphi^{k,i}\|_{L^2({\bf R}^N,e^{|x|^2/4}\,dx)}^2
\le C(1+t)^{-\frac{d_m}{2}}\|u_m(0)\|_{L^2({\bf R}^N,e^{|x|^2/4}\,dx)}^2
\end{split}
\end{equation*}
for $t>0$. This implies \eqref{eq:5.2}. 
On the other hand, 
by Lemma~\ref{Lemma:2.3} we have 
$$
\frac{|u_m(x,2t)|}{U(\min\{|x|,\sqrt{t}\})}\le Ct^{-\frac{N}{4}}U(\sqrt{t})^{-1}\|u_m(t)\|_{L^2({\bf R}^N)},
\quad x\in{\bf R}^N,\,\,t>0.
$$
This together with \eqref{eq:5.2} implies \eqref{eq:5.3}. 
Thus Lemma~\ref{Lemma:5.1} follows.
$\Box$
\vspace{5pt}
\newline
{\bf Proof of Theorem~\ref{Theorem:1.4}.}
Let $\varphi\in L^2({\bf R}^N,e^{|x|^2/4}\,dx)$ and $v:=e^{-tL_0}\varphi^{0,1}$. 
Let $K$ be any compact set in ${\bf R}^N\setminus\{0\}$ and $R>0$. 
In cases (S) and (C), by Theorems~\ref{Theorem:1.1} and \ref{Theorem:1.2}
we have 
\begin{equation}
\label{eq:5.4}
\begin{split}
 & \lim_{t\to\infty}t^{\frac{N+A}{2}}v\left(t^{\frac{1}{2}}y,t\right)=c_dm(\varphi^{0,1})|y|^Ae^{-\frac{|y|^2}{4}}
\,\,\,\mbox{in}\,\,\,L^2({\bf R}^N,e^{|y|^2/4}\,dy)\,\cap\,L^\infty(K),\\
 & \lim_{t\to\infty}\,t^{\frac{N+2A}{2}}\frac{v(x,t)}{U(|x|)}
 =\frac{c_d}{c_*}m(\varphi^{0,1})\,\,\,\mbox{in}\,\,\,L^\infty(B(0,R)).
\end{split}
\end{equation}
In case ($\mbox{S}_*$), Theorem~\ref{Theorem:1.3} implies that
\begin{equation}
\label{eq:5.5}
\begin{split}
 & \lim_{t\to\infty}t^{\frac{N+A}{2}}(\log t)v\left(t^{\frac{1}{2}}y,t\right)=2c_dm(\varphi^{0,1})|y|^Ae^{-\frac{|y|^2}{4}}
\,\,\,\mbox{in}\,\,\,L^2({\bf R}^N,e^{|y|^2/4}\,dy)\,\cap\,L^\infty(K),\\
 & \lim_{t\to\infty}\,t^{\frac{N+2A}{2}}(\log t)^2\frac{v(x,t)}{U(|x|)}
 =\frac{2\sqrt{2}}{c_*}m(\varphi^{0,1})
 =\frac{4c_d}{c_*}m(\varphi^{0,1})\,\,\,\mbox{in}\,\,\,L^\infty(B(0,R)).
\end{split}
\end{equation}
Here
\begin{equation}
\label{eq:5.6}
\begin{split}
c_dm(\varphi^{0,1}) & =\frac{c_d^2}{c_*}\int_0^\infty \varphi^{0,1}(r)U(r)r^{N-1}\,dr
=\frac{c_d^2}{c_*}|{\bf S}^{N-1}|^{-1}
\int_{{\bf R}^N} \varphi^{0,1}(|x|)U(|x|)\,dx\\
 & =\frac{c_d^2}{c_*}|{\bf S}^{N-1}|^{-1}
\int_{{\bf R}^N} \varphi(x)U(|x|)\,dx=M(\varphi). 
\end{split}
\end{equation}

Taking a sufficiently large integer $m$, 
by Lemma~\ref{Lemma:5.1} we have 
\begin{equation}
\label{eq:5.7}
\begin{split}
 & \lim_{t\to\infty}\,t^{\frac{N+A}{2}}u_m\left(t^{\frac{1}{2}}y,t\right)=0
\quad\mbox{in}\quad
L^2({\bf R}^N,e^{|y|^2/4}\,dy)\,\cap\,L^\infty(K),\\
 & \lim_{t\to\infty}\,t^{\frac{N+2A}{2}}\frac{u_m(x,t)}{U(|x|)}=0
 \qquad\,\mbox{in}\quad L^\infty(B(0,R)),
\end{split}
\end{equation}
for any compact set $K\subset{\bf R}^N\setminus\{0\}$ and $R>0$. 
On the other hand, 
$L_k$ is subcritical and $A^+(\lambda_2+\omega_k)>A$ for $k\in\{1,2,\dots,m-1\}$. 
Then, taking a sufficiently small $\epsilon>0$ if necessary, 
by Theorems~\ref{Theorem:1.1} and \ref{Theorem:1.2} we obtain 
\begin{equation}
\label{eq:5.8}
\begin{split}
 & \lim_{t\to\infty}\,t^{\frac{N+A}{2}}
\left[e^{-tL_k}\phi^{k,i}\right]\left(t^{\frac{1}{2}}y,t\right)=0\quad\mbox{in}\quad
L^2({\bf R}^N,e^{|y|^2/4}\,dy)\,\cap\,L^\infty(K),\\
 & \lim_{t\to\infty}t^{\frac{N+2A}{2}}\frac{\left[e^{-tL_k}\phi^{k,i}\right](x)}{U_k(|x|)}=0
\quad\,\,\qquad\mbox{in}\quad L^\infty(B(0,R)),
\end{split}
\end{equation}
for any compact set $K\subset{\bf R}^N\setminus\{0\}$ and $R>0$. 
On the other hand, it follows from \eqref{eq:1.2} that \eqref{eq:5.1} that 
$U_k(r)/U(r)$ is bounded on $(0,R)$ for any $R>0$. 
This together with \eqref{eq:5.8} implies that 
\begin{equation}
\label{eq:5.9}
\lim_{t\to\infty}\,t^{\frac{N+2A}{2}}\frac{\left[e^{-tL_k}\phi^{k,i}\right](|x|)}{U(|x|)}=0
\quad\mbox{in}\quad L^\infty(B(0,R))
\end{equation}
for any $R>0$. 
Since 
$$
[e^{-tL}\varphi](x)=v(x,t)+\sum_{k=1}^{m-1}\left[e^{-tL_k}\phi^{k,i}\right](|x|)Q_{k,i}\left(\frac{x}{|x|}\right)+u_m(x,t),
$$
by \eqref{eq:5.4}--\eqref{eq:5.9}
we obtain assertions~(a) and (b). 
Thus the proof is complete.
$\Box$
\vspace{5pt}
\newline
{\bf Proof of Corollary~\ref{Corollary:1.1}.}
Let $p=p(x,y,t)$ be the fundamental solution corresponding to $e^{-tL}$. 
Let $y\in{\bf R}^N$ and $\tau>0$. Set 
$\varphi(x)=p(x,y,\tau)$ for $x\in{\bf R}^N$. 
Taking a sufficiently small $\tau>0$ if necessary, 
by \eqref{eq:1.6} we see that 
$\varphi\in L^2({\bf R}^N,e^{|x|^2/4}\,dx)$. 
On the other hand, 
since $p(x,y,t)=p(y,x,t)$, 
we have
$$
\int_{{\bf R}^N}\varphi(x)U(|x|)\,dx
=\int_{{\bf R}^N}p(x,y,\tau)U(|x|)\,dx
=\int_{{\bf R}^N}p(y,x,\tau)U(|x|)\,dx
=U(|y|)
$$
for $y\in{\bf R}^N$ and $\tau>0$. 
Then, applying Theorem~\ref{Theorem:1.4} and letting $\tau\to+0$, 
we obtain the desired results. Thus Corollary~\ref{Corollary:1.1} follows. 
$\Box$
\vspace{5pt}

\noindent
{\bf Acknowledgements.} 
The first author was partially supported 
by the Grant-in-Aid for Scientific Research (A)(No.~15H02058)
from Japan Society for the Promotion of Science. 

\end{document}